\documentclass[12pt]{article}
\usepackage{amsfonts,amsmath,amsthm}
\usepackage{amssymb,latexsym}
\usepackage[dvips]{epsfig}
\usepackage{psfrag}
\usepackage{amscd}
\usepackage{pstricks,pst-node}
\usepackage{rotating}

\title{Categorification of level two representations of quantum $sl_n$ via generalized arc rings}

\author{Yanfeng Chen}
\date{\today}

\newtheorem{prop}{Proposition}
\newtheorem{theorem}{Theorem}
\newtheorem{lemma}{Lemma}

\newcommand{\oplusop}[1]{{\mathop{\oplus}\limits_{#1}}}
\newcommand{\mo}{\mathbf{1}}
\newcommand{\cA}{{\mathcal{A}}}
\newcommand{\cF}{{\mathcal{F}}}
\newcommand{\cK}{{\mathcal{K}}}
\newcommand{\Inv}{\textrm{Inv}}
\newcommand{\define}{\stackrel{\mbox{\scriptsize{def}}}{=}}
\def\sbinom#1#2{\left( \hspace{-0.06in}\begin{array}{c} #1 \\ #2 \end{array}
 \hspace{-0.06in} \right)}
 \def\leftrightmaps#1#2#3{\raise3pt\hbox{$\mathop{\,\,\hbox to
     #1pt{\rightarrowfill}\kern-#1pt\lower3.95pt\hbox to
     #1pt{\leftarrowfill}\,\,}\limits_{#2}^{#3}$}}
\def\mc{\mathcal}

\def\cA{\mc{A}}
\def\cC{\mc{C}}
\def\cE{\mc{E}}
\def\cF{\mc{F}}
\def\cK{\mc{K}}

\def\id{\mathrm{Id}}

\def\hsp{\hspace{0.1in}}

\begin{document}

\maketitle \baselineskip 14pt \setcounter{tocdepth}{1}
\tableofcontents

\section{Introduction}
Khovanov constructed in \cite{Khovanov:FunctorValuedInvariant} a
family of rings $H^n$, for $n\geq 0$, which is a categorification
of \Inv($n$), the space of $U_q(sl_2)$-invariants in $V^{\otimes
2n}.$ These rings lead to an invariant of (even) tangles which to
a tangle assigns a complex of $(H^n,H^m)$-bimodules, up to chain
homotopy equivalence. Khovanov and the author
\cite{CK:Subquotients} built subquotients of $H^n$ and used them
to categorify the action of tangles on $V^{\otimes n}$. The same
rings were also introduced by Stroppel \cite{S:PerverseSheaves}.

In this paper, we extend the construction of these arc rings
$A^{k,n-k}$ and give a categorification of level two
representations of $U_q(sl_N)$. In section $2$ we review the
definition of the arc rings $A^{n-k,k}$ and construct the rings
$A_n^{k,l}$ with two platforms of arbitrary sizes $k$ and $l$. We
show that they lead to a tangle invariant which is functorial
under tangle cobordisms. Then in section $3$ we compute the
centers of $A_n^{k,l}$ and relate them to the cohomology rings of
Springer varieties. Finally, in section $4$, we categorify level
two representations of $U_q(sl_N)$ using the rings $A_n^{k,0}$.

Fix a level two representation $V$ of $U_q(sl_N)$ with the highest
weight $\omega_{s}+\omega_{k+s}$. There is a decomposition of $V$
into weight spaces $V= \oplusop{\mu} V_{\mu}$. A weight $\mu$ is
called admissible if it appears in the weight space decomposition
of $V$. Denote by $\mathcal{C}$ the direct sum of categories of
$A^{k,0}_{m(\mu)}$-modules over admissible $\mu$, where $m(\mu)$
is a non-negative integer depending only on $\mu$. The
Grothendieck group of the category of $A^{k,0}_{m(\mu)}$-modules
is naturally isomorphic to $V_{\mu}$. The exact functors
$\mathcal{E}_i,\ \mathcal{F}_i$ introduced by Khovanov and
Huerfano in \cite{HK:LevelTwo} naturally extend to exact functors
on $\mathcal{C}$ which categorify the actions of $E_i,F_i\in
U_q(sl_N)$ on $V$.

{\bf Acknowledgements:} I would like to thank Mikhail Khovanov for
recommending this problem and too many helpful conversations and
suggestions. I am also very grateful to Robin Kirby for his
kindness, guidance and support.

\section{Generalization of the arc ring $A^{n-k,k}$}
\subsection{Arc ring $A^{n-k,k}$}
We first recall the definition of $H^n$ from
\cite{Khovanov:FunctorValuedInvariant}. Let $\cA$ be a free
abelian group of rank $2$ spanned by $\mo$ and $X$ with $\mo$ in
degree $-1$ and $X$ in degree $1$. Assign to $\cA$ a
$2$-dimensional TQFT $\cF$ which associates $\cA^{\otimes k}$ to a
disjoint union of $k$ circles. To the ``pants'' cobordism
corresponding to merging of two circles into one, $\cF$ associates
the multiplication $m: \cA \otimes \cA \rightarrow \cA$
 \begin{equation}
 \mo^2=\mo, \hspace{0.1in} \mo X= X\mo = X, \hspace{0.1in} X^2=0. \label{m}
 \end{equation}
To the ``inverse pants'' cobordism corresponding to splitting of
one circle into two, $\cF$ associates the comultiplication $
\Delta: \cA \rightarrow \cA \otimes\cA$
 \begin{equation}
\hspace{0.1in}\Delta(\mo) = \mo \otimes X + X\otimes \mo ,
\hspace{0.1in} \Delta(X) = X\otimes
 X.\label{Delta}
 \end{equation}
 To the ``cup'' and ``cap'' cobordisms corresponding to the birth and death of a circle, $\cF$ associates the unit map $\iota: \mathbb{Z}
\rightarrow \cA$ and trace map $\epsilon: \cA \rightarrow
\mathbb{Z}$ respectively
 \begin{equation*}
 \iota(1)=\mo,\hspace{0.1in} \varepsilon(\mo)=0, \hspace{0.1in}\varepsilon(X)=1.
 \end{equation*}

Let $B^n$ be the set of crossingless matchings of $2n$ points. For
$a$, $b\in B^n$ denote by $W(b)$ the reflection of $b$ about the
horizontal axis, and by $W(b)a$ the closed $1$-manifold obtained
by closing $W(b)$ and $a$ along their boundaries.

\begin{figure}[ht!]
\begin{center}
\psfrag{b}{\small$b$}\psfrag{a}{\small$a$}\psfrag{W(b)}{\small$W(b)$}\psfrag{W(b)a}{\small$W(b)a$}
\psfrag{f}{\small$\cF$}\psfrag{a2}{\small$\cA^{\otimes 2}$}
\epsfig{figure=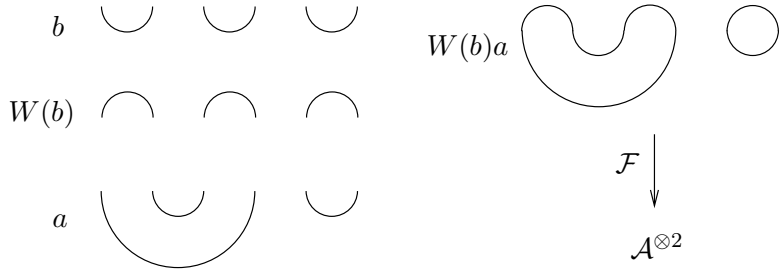}  \caption{Gluing in $B^3$.}
\label{Wba.figure}
\end{center}
\end{figure}

$\cF(W(b)a)$ is a graded abelian group isomorphic to $\cA^{\otimes
I}$, where $I$ is the set of circles in $W(b)a$, see
figure~\ref{Wba.figure}.

\begin{figure}[ht!]
\begin{center}
\psfrag{a}{\tiny$a$}\psfrag{c}{\tiny$c$}\psfrag{W(b)}{\tiny$W(b)$}\psfrag{W(a)}{\tiny$W(a)$}
\psfrag{f}{\tiny$\cF$}\psfrag{a1}{\tiny$\cA$}\psfrag{a4}{\tiny$\cA^{\otimes
4}$}\psfrag{m}{\tiny$m\circ(id\otimes m)\circ(id^2\otimes m)$}
\epsfig{figure=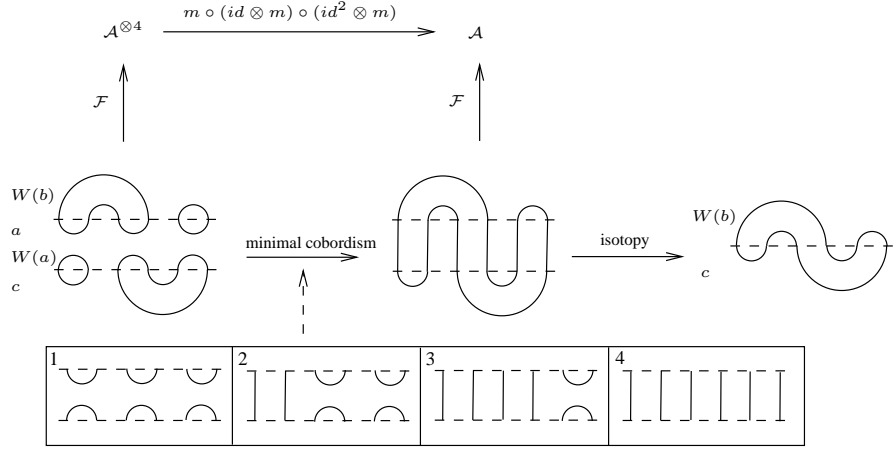} \caption{Multiplication in $H^n$.}
\label{contraction.figure}
\end{center}
\end{figure}

For $a$, $b\in B^n$ let
 \begin{equation*}
 {_b(H^n)_a} \define \cF(W(b)a)\{n\}.
 \end{equation*}
and define $H^n$ as the direct sum
 \begin{equation*}
 H^n\define \oplusop{a,b} \hspace{0.05in} {_b(H^n)_a},
 \end{equation*}
where $\{n\}$ denotes the action of raising the grading up by $n$.
Multiplication maps in $H^n$ are defined as follows. We set $xy=0$
if $x\in {_b(H^n)_a}$, $y\in {_c(H^n)_d}$ and $c\neq a$.
Multiplication maps
 \begin{equation*}
 {_b(H^n)_a} \otimes {_a(H^n)_c} \rightarrow {_b(H^n)_c}
 \end{equation*}
are given by homomorphisms of abelian groups
 \begin{equation*}
 \cF(W(b)a) \otimes \cF(W(a)c) \rightarrow \cF(W(b)c),
 \end{equation*}
which are induced by the ``minimal'' cobordism from $W(b)aW(a)c$
to $W(b)c$, see figure~\ref{contraction.figure}.

\begin{figure}[ht!]
\begin{center}
\epsfig{figure=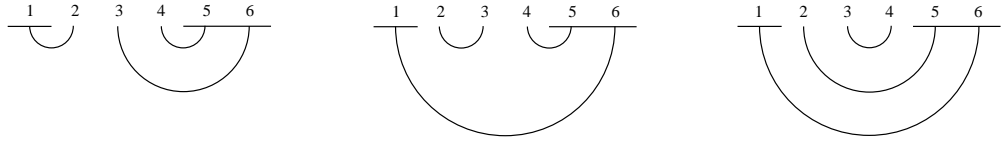}  \caption{The $3$ elements in $B^{1,2}$.}
\label{b21.figure}
\end{center}
\end{figure}

Now we recall the definition of the subquotients of $H^n$ from
\cite{CK:Subquotients}. For each $n\geq 0$ and $0\leq k \leq n$,
define $B^{n-k,k}$ to be the subset of $B^n$ where there are no
matchings among the first $n-k$ points and among the last $k$
points. Figure~\ref{b21.figure} shows $B^{1,2}$.
We put two ``platforms'', one on the first $n-k$ points and one on
the last $k$ points to indicate that these endpoints are special.
The $n$ points lying in between the two platforms are called ``free points''.

Define $\widetilde{A}^{n-k,k}$ by
 \begin{equation}
   \widetilde{A}^{n-k,k} \define \oplusop{a,b\in B^{n-k,k}}
   \hspace{0.05in} \cF(W(b)a)\{n\}. \label{DefTildeA.equation}
 \end{equation}
$\widetilde{A}^{n-k,k}$ sits inside $H^n$ as a graded subring
which inherits its multiplication from $H^n$.

\begin{figure}[ht!]
\begin{center}
\psfrag{1}{\tiny I} \psfrag{2}{\tiny II} \psfrag{3}{\tiny III}
\epsfig{figure=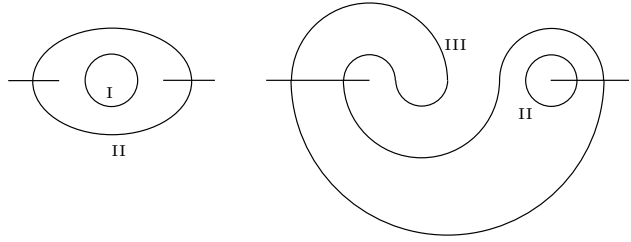} \caption{$3$ types of circles.}
\label{3types.figure}
\end{center}
\end{figure}

For $a,b\in B^{n-k,k}$, call a circle in $W(b)a$ type I if it is
disjoint from platforms, type II if it intersects at least one
platform and intersect each platform at most once, and type III if
it intersects one of the platforms at least twice (see
figure~\ref{3types.figure}). An intersection point between a
circle and a platform is called a ``mark''.

Ring $\widetilde{A}^{n-k,k}$ has a two-sided graded ideal
$I^{n-k,k}\subset \widetilde{A}^{n-k,k}$. Ring $A^{n-k,k}$ is
defined as the quotient of $\widetilde{A}^{n-k,k}$ by the ideal
$I^{n-k,k}$
 \begin{equation}
   A^{n-k,k} \define \widetilde{A}^{n-k,k} / I^{n-k,k}. \label{DefA.equation}
 \end{equation}
 $A^{n-k,k}$ naturally decomposes into a direct sum of graded
 abelian groups
 \begin{equation*}
   A^{n-k,k} =  \oplusop{a,b\in B^{n-k,k}} \hspace{0.05in}
   _a(A^{n-k,k})_b.
 \end{equation*}
By taking the direct product over all $0 \leq k \leq n$, we
collect the rings $A^{n-k,k}$ together into a graded ring $A^{n}$
$$A^n\define \prod_{0\leq k\leq n}A^{n-k,k}.$$
As a graded abelian group, $A^n$ is the direct sum of $A^{n-k,k}$,
over $0\leq k\leq n$.

See \cite{Khovanov:FunctorValuedInvariant} and
\cite{CK:Subquotients} for more details on $H^n$ and its
subquotients.

\subsection{Generalization of $A^{n-k,k}$}
We call the triple $(n,k,l)$ coherent $|k-l|\leq n$ and
$n+k+l\equiv 0\ \ \mathrm{(mod\ 2)}$. For each coherent triple
$(n,k,l)$ define $B_n^{k,l}$ to be the subset of $B^{(n+k+l)/2}$
where there are no matchings among the first $k$ points and among
the last $l$ points. Put one platform on the first $k$ points and
one on the last $l$ points. Note that $B_{2n}^{0,0}=B^n$ and
$B_n^{n-k,k}=B^{n-k,k}$.

Define $\widetilde{A}_n^{k,l}$ by
 \begin{equation}
   \widetilde{A}_n^{k,l} \define \oplusop{a,b\in B_n^{k,l}}
   \hspace{0.05in} \cF(W(b)a)\{\frac{n+k+l}{2}\}. \label{DefTildeAnkl.equation}
 \end{equation}
Just like $\widetilde{A}^{n-k,k}$, $\widetilde{A}_n^{k,l}$ is a
graded subring of $H^n$ and inherits its multiplication from
$H^n$. The ideal $I_n^{k,l}\subset \widetilde{A}_n^{k,l}$ is
defined exactly as that of $\widetilde{A}^{n-k,k}$. For $a,b\in
B_n^{k,l}$, if $W(b)a$ contains at least one type III circle, set
${_b(I_n^{k,l})_a}=\cF(W(b)a)\{n\}$. If $W(b)a$ contains only
circles of type I and type II, we write $\cF(W(b)a)=\cA^{\otimes
i}\otimes \cA^{\otimes j}$ in which type II circles correspond to
the first $i$ tensor factors, and define ${_b(I_n^{k,l})_a}$ as
the span of
 \begin{equation*}
y_1 \otimes \cdots \otimes y_{t-1} \otimes X \otimes y_{t+1}
\otimes \cdots \otimes y_{i+j} \in   \cA^{\otimes i}\otimes
\cA^{\otimes j} \cong \cF(W(b)a),
 \end{equation*}
where $1\leq t\leq i$ and $y_s\in \{\mo, X\}$. By taking the
direct sum over all $a,b\in B_n^{k,l}$ we get a subgroup of
$\widetilde{A}_n^{k,l}$
 \begin{equation*}
   I_n^{k,l} \define \oplusop{a,b\in B_n^{k,l}}\hspace{0.05in}
   {_b(I_n^{k,l})_a}.
 \end{equation*}
It's easy to show that $I_n^{k,l}$ is a two-sided graded ideal of
$\widetilde{A}_n^{k,l}$. Ring $A_n^{k,l}$ is defined as the
quotient of $\widetilde{A}_n^{k,l}$ by the ideal $I_n^{k,l}$
 \begin{equation}
   A_n^{k,l} \define \widetilde{A}_n^{k,l} / I_n^{k,l}. \label{DefAnkl.equation}
 \end{equation}
If $W(b)a$ contains a type III circle then ${_a(A^{n-k,k})_b}=0$.
Otherwise, group ${_a(A^{n-k,k})_b}$ is free abelian of rank
$2^{\mathrm{\#\ of\ type\ I\ circles}}$. Assuming that $W(a)b$
contains $m$ circles in which the first $i$ of them are of type
II, ${_a(A_n^{k,l})_b}$ has a basis of the form
 \begin{equation*}
\mo \otimes \cdots \otimes \mo \otimes a_{i+1} \otimes \cdots
\otimes a_{m},
 \end{equation*}
 where $a_s\in \{\mo, X\}$ for all $i+1\leq s\leq m$.

There is a natural decomposition of $A_n^{k,l}$ into a direct sum
of graded abelian groups
 \begin{equation*}
   A_n^{k,l} =  \oplusop{a,b\in B_n^{k,l}} \hspace{0.05in}
   _a(A_n^{k,l})_b,
 \end{equation*}
where $$_a(A_n^{k,l})_b = \cF(W(a)b)\{\frac{n+k+l}{2}\} /
{_a(I_n^{k,l})_b}.$$ Let $P_n^{k,l}(a)$, or simply P(a), for $a\in
B_n^{k,l}$, be a left $A_n^{k,l}$-module given by
 $$P(a)=\oplusop{b\in B_n^{k,l}} {_b(A^{n-k,k})_a}.$$
$A_n^{k,l}$ decomposes into a direct sum of left
$A_n^{k,l}$-modules
$$A_n^{k,l}=\oplusop{a\in B_n^{k,l}}P(a).$$
$P(a)$ is left projective since it is a direct summand of the free
module $A_n^{k,l}$. Actually, any indecomposable left projective
$A_n^{k,l}$-module is isomorphic to $P(a)\{s\}$ for some $a\in
B_n^{k,l}$ and $s\in\mathbb{Z}$.

Here are some basic facts about the ring $A_n^{k,l}$:
\begin{itemize}
\item $A_n^{k,l}\cong A_n^{l,k}$. Reflecting a diagram in
$B_n^{k,l}$ about a vertical axis produces a diagram in
$B_n^{l,k}$. It leads to an isomorphism of sets $B_n^{k,l}\cong
B_n^{l,k}$ which induces an isomorphism of rings
$\widetilde{A}_n^{k,l}\cong \widetilde{A}_n^{l,k}$ and of the
quotient rings $A_n^{k,l}\cong A_n^{l,k}$.

\item The minimal idempotents in $A_n^{k,l}$ are $1_a\define
\mo^{\otimes (n+k+l)/2}\in {_a(A_n^{k,l})_a}$. The unit element
$1$ of $A_n^{k,l}$ is the sum of $1_a$ over all $a \in B_n^{k,l}$:
$1 \define \sum_{a\in B_n^{k,l}} 1_a.$

\item $A_n^{k,l}$ sits inside $A_{n}^{k+1,l+1}$ as a subring. This
inclusion stabilizes when $k+l>n$. In particular we have:
$A_n^{k,n-k}\cong A_n^{k+1,n-k+1}\cong A_n^{k+2,n-k+2}\cong
\cdots$
\end{itemize}

\begin{prop}
The rings $A_n^{0,l}$ are symmetric and, therefore, Frobenius over
$\mathbb{Z}$.
\end{prop}
The proof is similar to \cite{Khovanov:FunctorValuedInvariant}
proposition $32$.

\subsection{Flat tangles and bimodules}

\begin{figure}[ht!]
\begin{center}
\epsfig{figure=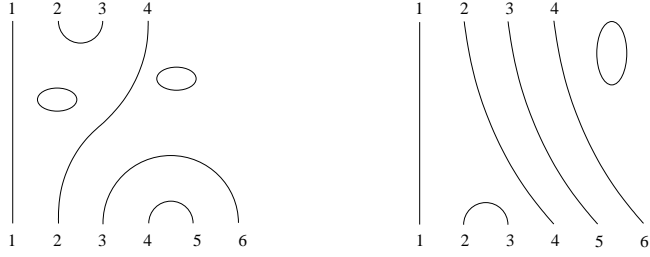}  \caption{Two flat tangles
in $\widehat{B}^4_6$.} \label{FlattangleExample.figure}
\end{center}
\end{figure}

Denote by $\widehat{B}^m_n$ the space of flat tangles with $m$ top
endpoints and $n$ bottom endpoints. For simplicity we assume that
the top and bottom endpoints lie on $\mathbb{R}\times \{1\}$ and
$\mathbb{R}\times \{0\}$, and have integer coefficients
$1,2,\cdots,m$ and $1,2,\cdots,n$ respectively.
Figure~\ref{FlattangleExample.figure} shows two elements in
$\widehat{B}^4_6$.

\begin{figure}[ht!]
\begin{center}
\psfrag{a}{\small $c\in B_2^{0,2}$} \psfrag{t}{\small $T\in
\widehat{B}_6^2$} \psfrag{b}{\small $b\in
B_6^{1,3}$}\psfrag{atb}{\small $W(c)Tb$} \psfrag{c}{\small closure
of $W(c)Tb$} \psfrag{v}{\small Vertical lines added}
 \epsfig{figure=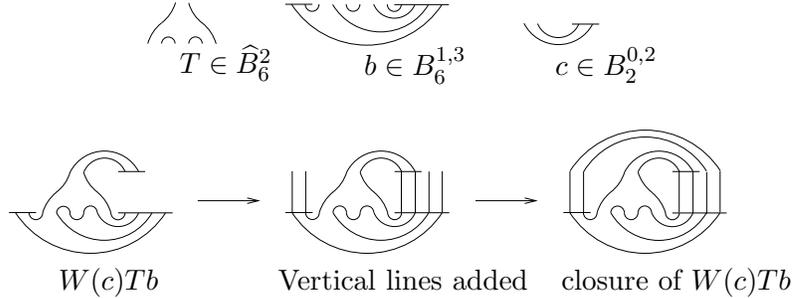}  \caption{Closing $W(c)Tb$.} \label{close.figure}
\end{center}
\end{figure}

To a flat tangle $T\in \widehat{B}^m_n$ we would like to assign a
bimodule over algebras $A_n^{k,l}$ and $A_m^{s,t}$ where both
$(n,k,l)$ and $(m,s,t)$ are coherent triples and $k-l=s-t$. Define
a graded $(\widetilde{A}_m^{s,t},\widetilde{A}_n^{k,l})$-bimodule
$\widetilde{\cF}(T)$ by
  \begin{equation*}
    \widetilde{\cF}(T) = \oplusop{b\in B_n^{k,l},c\in B_m^{s,t}} \hspace{0.05in}{_c\widetilde{\cF}(T)_b},
  \end{equation*}
where
  \begin{equation*} \label{def-bimod}
    {_c}\widetilde{\cF}(T)_b \stackrel{\mbox{\scriptsize{def}}}{=} \cF( W(c) T b)
    \{\frac{n+k+l}{2}\}.
  \end{equation*}
The plane diagram $W(c) T b$ is not a union of circles if $k\neq
s$. In that case we close it in the obvious way before applying
the functor $\cF$ (see figure~\ref{close.figure}). The left action
$\widetilde{A}_m^{s,t} \times \widetilde{\cF}(T) \rightarrow
\widetilde{\cF}(T)$ comes from maps
  \begin{equation*}
    \cF(W(a)c) \times {_c\widetilde{\cF}(T)_b} \rightarrow
    {_a\widetilde{\cF}(T)_b},
  \end{equation*}
and the right action $\widetilde{\cF}(T) \times
\widetilde{A}_n^{k,l} \rightarrow \widetilde{\cF}(T)$ comes from
maps
  \begin{equation*}
    {_c\widetilde{\cF}(T)_b} \times \cF(W(b)a)\rightarrow
    {_c\widetilde{\cF}(T)_a}.
  \end{equation*}
Both maps are induced by the obvious ``minimal cobordism'' (see
figure~\ref{contraction.figure}).

Now define a subgroup ${_bI(T)_a}$ of ${_b\widetilde{\cF}(T)_a}$
as follows. Set ${_bI(T)_a}={_b\widetilde{\cF}(T)_a}$ if $W(b)Ta$
contains a type III arc. Otherwise, assuming that
$\cF(W(b)Ta)\cong \cA^{\otimes r}$ in which type II circles
correspond to the first $i$ tensor factors, set ${_bI(T)_a}$ to be
the span of
 \begin{equation*}
u_1 \otimes \cdots \otimes a_{j-1} \otimes X \otimes u_{j+1}
\otimes \cdots \otimes u_{r} \in \cF(W(b)Ta) \cong \cA^{\otimes
r},
 \end{equation*}
  where $1\leq j\leq i$ and $u_e\in \{\mo, X\}$ for each $1\leq e\leq r$, $e\neq j$. By taking the direct sum we get a subgroup
 \begin{equation*}
   I(T) \define \oplusop{a \in B_m^{s,t}, b\in B_n^{k,l}}\hspace{0.05in}
   {_aI(T)_b}.
 \end{equation*}
$I(T)$ is in fact a subbimodule of $\widetilde{\cF}(T)$ and we can
define $\cF(T)$ to be the quotient bimodule
 \begin{equation*}
   \cF(T) \define \widetilde{\cF}(T) / I(T).
 \end{equation*}
It's easy to show that the action of $I_n^{k,l}$ on $\cF(T)$ is
trivial (see \cite{CK:Subquotients}), thus the
$(\widetilde{A}_m^{s,t},\widetilde{A}_n^{k,l})$-bimodule structure
on $\cF(T)$ descends to an $(A_m^{s,t},A_n^{k,l})$-bimodule
structure.

\begin{prop} An isotopy between $T_1,T_2\in \widehat{B}^{m}_{n}$ induces an isomorphism
 of bimodules $\cF(T_1) \cong \cF(T_2).$ Two isotopies between $T_1$ and $T_2$
 induce equal isomorphisms iff the bijections from circle components
 of $T_1$ to circle components of $T_2$ induced by the two isotopies coincide.
\end{prop}
The proof is similar to that in
\cite{Khovanov:FunctorValuedInvariant}.

\begin{figure}[ht!]
\begin{center}
\psfrag{c}{\small cobordism}\psfrag{a1}{\small
$a$}\psfrag{b1}{\small $b_1$}\psfrag{b2}{\small
$b_2$}\psfrag{b3}{\small $b_3$}\psfrag{t1}{\small
$T_1$}\psfrag{t2}{\small $T_2$}\psfrag{e1}{\tiny $1\otimes 1
\otimes \mo \mapsto  1\otimes 1$}\psfrag{e2}{\tiny $1\otimes 1
\otimes X\mapsto 0$}\psfrag{e3}{\tiny $1\otimes 1\mapsto X \otimes
1\otimes 1$}\psfrag{e4}{\tiny $0\mapsto 0$}\psfrag{w1}{\tiny
$W(a)T_1 b_1$}\psfrag{w2}{\tiny $W(a)T_1 b_2$}\psfrag{w3}{\tiny
$W(a)T_1 b_3$}\psfrag{w4}{\tiny $W(a)T_2 b_1$}\psfrag{w5}{\tiny
$W(a)T_2 b_2$}\psfrag{w6}{\tiny $W(a)T_2 b_3$}\psfrag{f}{\small
$\cF$}\psfrag{a}{\small $\cA \{1\}$}\psfrag{z}{\small $\mathbb{Z}
\{1\}$}\psfrag{0}{\small $0$}\psfrag{1}{\tiny $1$}
\psfrag{x}{\small $1$}\psfrag{y}{\small $0$}
 \epsfig{figure=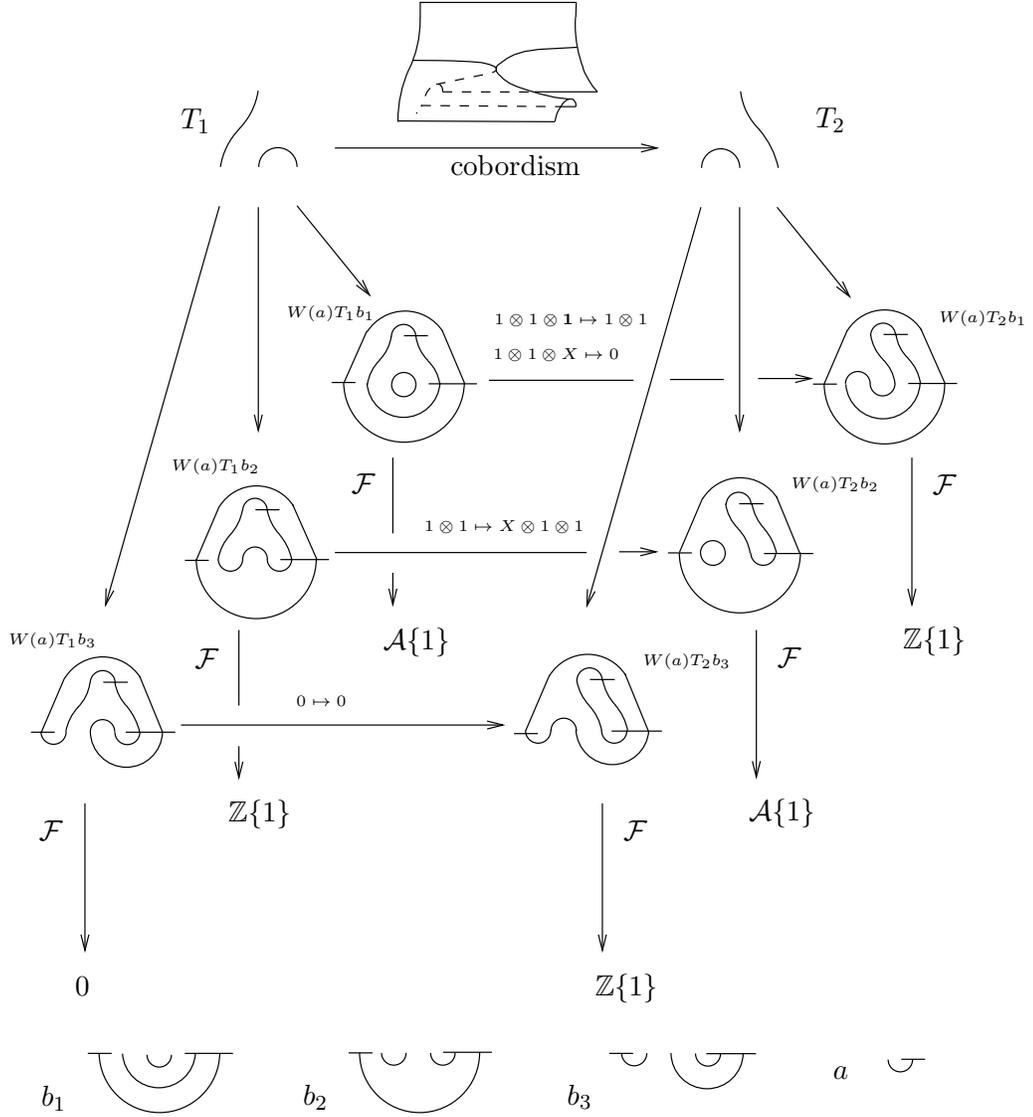}  \caption{Cobordism induces bimodule map} \label{Cobordism.figure}
\end{center}
\end{figure}

Cobordisms between flat tangles induce bimodule maps (see
figure~\ref{Cobordism.figure}):
\begin{prop} Let $T_1,T_2\in \widehat{B}^{m}_{n}$ and $S$ a cobordism between $T_1$ and $T_2$.
Then $S$ induces a degree $\frac{n+m}{2}-\chi(S)$ homomorphism of
$(A_m^{s,t},A_n^{k,l})$-bimodules
   \begin{equation*}
     \cF(S): \cF(T_1) \to \cF(T_2),
   \end{equation*}
where $\chi(S)$ is the Euler characteristic of $S$.
\end{prop}
\begin{proof} It follows from the definition that
$\widetilde{\cF}(T_1)= \oplusop{a,b} \cF(W(b)T_1 a) \{
\frac{n+k+l}{2}\}$ and $\widetilde{\cF}(T_2)= \oplusop{a,b}
\cF(W(b)T_2 a) \{ \frac{n+k+l}{2}\}$, where the sum is over all
$a\in B_n^{k,l}$ and $b\in B_m^{s,t}$. The surface $S$ induces a
homogeneous map of graded abelian groups $\cF(W(b)T_1 a) \to
\cF(W(b)T_2 a)$. Summing over all $a$ and $b$ we get a
homomorphism of
$(\widetilde{A}_n^{k,l},\widetilde{A}_m^{s,t})$-bimodules
$\widetilde{\cF}(S): \widetilde{\cF}(T_1) \to
     \widetilde{\cF}(T_2)$. The grading assertion follows from the fact that $\chi(S')=\chi(S)-\frac{n+m}{2}$.
It's easy to show that $\widetilde{\cF}(S)$ takes $I(T_1)$ into
$I(T_2)$. See \cite{CK:Subquotients} for details.
\end{proof}

\begin{prop} Isotopic (rel boundary) surfaces induce equal bimodule maps.
\end{prop}

\begin{prop} Let $T_1,T_2,T_3\in \widehat{B}^{m}_{n}$ and $S_1$,
$S_2$ be cobordisms from $T_1$ to $T_2$ and from $T_2$ to $T_3$
respectively. Then $\cF(S_2)\cF(S_1)=\cF(S_2\circ S_1)$.
\end{prop}

Proofs of the above two propositions are similar to those in
\cite{Khovanov:FunctorValuedInvariant}.

Two coherent triples $(n,k,l)$ and $(m,s,t)$ are called
\emph{compatible} if either $k+l=n$, $s+t=m$, $t=l+\frac{m-n}{2}$,
or $k=s$, $l=t$. Also, a $(A_m^{s,t},A_n^{k,l})$-bimodule is
called compatible if $(n,k,l)$ and $(m,s,t)$ are compatible.

\begin{prop} \label{projectivity.proposition}
Let $T\in\widehat{B}_n^m$, bimodule $\cF(T)$ is projective as a
left $A_m^{s,t}$-module and as a right $A_n^{k,l}$-module if
$(n,k,l)$ and $(m,s,t)$ are compatible.
\end{prop}
\begin{proof} Ignore all the grading shifts. The bimodule $\cF(T)$
is isomorphic, as a left $A_m^{s,t}$-module, to the direct sum
$\oplus_{a\in B_n^{k,l}} \cF(Ta)$. To prove $\cF(T)$ is left
projective it suffices to prove that $\cF(Ta)$ is left projective
for all $a\in B_n^{k,l}$. Fix any $a\in B_n^{k,l}$. In general,
$Ta$ is a union of circles and arcs. With all circles removed,
$Ta$ is isotopic to some $a'\in B_m^{k,l}$ (see
figure~\ref{deform.figure}).
\begin{figure}[ht!]
\begin{center}
\psfrag{a}{\small $a\in B_6^{3,1}$} \psfrag{T}{\small $T\in
\widehat{B}_6^4$} \psfrag{Ta}{\small $Ta$}\psfrag{a'}{\small
$a'\in B_4^{3,1}$} \psfrag{rm}{\small
Isotopy}\psfrag{iso}{$\cong$}
 \epsfig{figure=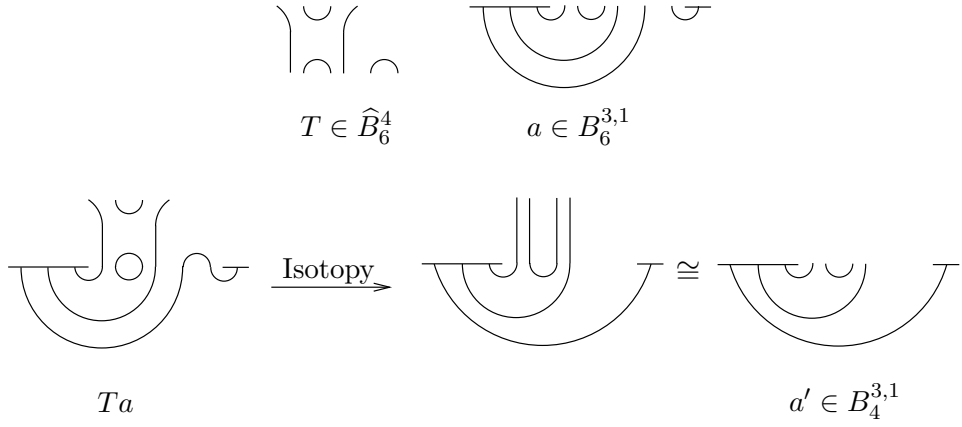}  \caption{Deformation of $Ta$.} \label{deform.figure}
\end{center}
\end{figure}

Case $1$: $k=s$ and $l=t$. In this case, assuming there are $c$ circles in $Ta$,
$$\cF(Ta)=(\oplusop{b\in A_m^{s,t}} \cF(W(b)Ta))\otimes {\cA}^{\otimes c}\cong (\oplusop{b\in A_m^{s,t}}\cF(W(b)a'))\otimes {\cA}^{\otimes c}.$$
By definition $\oplusop{b\in A_m^{s,t}}\cF(W(b)a')=
P_m^{s,t}(a')$, therefore $\cF(Ta)$ is left projective.

Case $2$: $k+l=n$, $s+t=m$, $t=l+\frac{m-n}{2}$. Without loss of
generality we assume that $m\leq n$. Let $\Theta=\frac{n-m}{2}$.
The case $\Theta=0$ is proved in case $1$. Suppose the statement
is true when $\Theta\leq d$. Consider any $T\in \widehat{B}_n^m$
and $a\in B_n^{k,l}$ such that $\frac{n-m}{2}=d+1$. There exists
at least one ``cap'' in $T$ which connects two bottom endpoints
since $n>m$. Pick a cap $c$ which has no other bottom endpoints of
$T$ between its two feet. After gluing $a$ to $T$, either there is
an arc in $a$ connecting the two platforms, or both feet of $c$ is
connected to the platforms since $k+l=n$. Therefore there is
always an arc connecting the two platforms in $Ta$. By definition
of the arc ring, the two far ends of the two platforms are then
connected by an arc $e$. When closing the graph $W(b)Ta$ for some
$b\in B_m^{s,t}$ we need to add $d+1$ arcs since $t<l$ (see
figure~\ref{close.figure}). Denote the topmost added arc by $f$.
The arcs $e$ and $f$ form a type II circle $g$ which encloses the
rest of $W(b)Ta$. We can remove $g$ from $W(b)Ta$ for any $b\in
B_m^{s,t}$ since it contributes nothing to $\cF(W(b)Ta)$, and then
reduce to the case $\Theta=d$. The proposition follows by
induction.
\end{proof}

\begin{prop} Let $T_1\in \widehat{B}^{p}_{n}$, $T_2\in
\widehat{B}^{m}_{p}$, $\cF(T_1)$ a compatible
$(A_{p}^{q,r},A_{n}^{k,l})$-bimodule, and $\cF(T_2)$ a compatible
$(A_{m}^{s,t},A_{p}^{q,r})$-bimodule. Then there is a canonical
isomorphism of $(A_{m}^{s,t},A_{n}^{k,l})$-bimodules
   \begin{equation*}
     \cF(T_2 T_1)\cong \cF(T_2) \otimes_{A_p^{q,r}}\cF(T_1).
   \end{equation*}
\end{prop}
\begin{proof} It follows from proposition
\ref{projectivity.proposition} that $W(a)T_2$ is a projective
right $A_p^{q,r}$-module and $T_1 b$ is a projective left
$A_p^{q,r}$-module for $a\in B_m^{s,t}$ and $b\in B_n^{k,l}$. The
proof in \cite{Khovanov:FunctorValuedInvariant} theorem $1$
therefore works in our case without any changes.
\end{proof}

{\bf{Remark:}} Compatible triples fall into two different types. A
pair of coherent triples $(n,k,l)$ and $(m,s,t)$ are called
``T-compatible'' if $k+l=n$, $s+t=m$, $t=l+\frac{m-n}{2}$, and
``F-compatible'' if $k=s$, $l=t$. Similarly, call a
$(A_n^{k,l},A_m^{s,t})$-bimodule T-compatible (F-compatible) if
$(n,k,l)$ and $(m,s,t)$ are T-compatible (F-compatible). If two
flat tangles $T_1$ and $T_2$ belong to the same type, $T_2 T_1$ is
then compatible and also belongs to that type. Therefore we can
compose as many flat tangles as we want within the same type.
However, if $T_1$ and $T_2$ belong to different types their
composition $T_2T_1$ may not be compatible.

Now consider only F-compatible triples and bimodules. For each $n$
such that $(n,k,l)$ is coherent, denote by $A_n^{k,l}$-mod the
category of finitely-generated graded left $A_n^{k,l}$-modules and
module maps. For each $T\in \widehat{B}_n^m$, tensoring with the
$(A_m^{k,l},A_n^{k,l})$-bimodule $\cF(T)$ is an exact functor from
$A_n^{k,l}$-mod to $A_m^{k,l}$-mod. A cobordism $S$ between two
flat tangles $T_1,T_2\in \widehat{B}_n^m$ induces a homomorphism
$\cF(S)$ of $(A_m^{k,l},A_n^{k,l})$-bimodules. The following
proposition is a summary of this section.

\begin{prop} Bimodules $\cF(T)$ and homomorphisms $\cF(S)$
assemble into a $2$-functor from the $2$-category of flat tangle
cobordisms to the $2$-category of natural transformations between
exact functors between $A_n^{k,l}$-mod.
\end{prop}

\subsection{Tangles and complexes of bimodules}

An $(m,n)$-tangle $L$ is a proper embedding of
$\frac{n+m}{2}$ oriented arcs and a finite number of oriented circles into
$\mathbb{R}^2 \times [0,1]$ such that the boundary points of arcs map to
\begin{equation*}
\{1, 2, . . ., n\}\times \{0\} \times\{0\}, \{1, 2, . . .,
m\}\times \{0\} \times\{1\}.
\end{equation*}
A plane diagram of a tangle is a generic projection of the tangle
onto $\mathbb{R}\times [0,1]$.

Fix $k$ and $l$ throughout the rest of this section. We would like
to define a tangle invariant using the rings $A_n^{k,l}$. The
construction follows the same line as in \cite{CK:Subquotients}.
The sizes of the platforms don't matter. We will state the results
here for completeness and refer readers to \cite{CK:Subquotients}
and \cite{Khovanov:FunctorValuedInvariant} for details.

Fix a diagram $D$ of an  oriented $(m,n)$-tangle $L$. We define
the complex of $(A_m^{k,l}, A_n^{k,l})$-bimodules $\cF(D)$
associated to $D$ inductively  as follows.

\begin{figure}[ht!]
\begin{center}
\psfrag{n}{Negative} \psfrag{p}{Positive} \psfrag{0}{\tiny
$0$-smoothing} \psfrag{1}{\tiny $1$-smoothing} \psfrag{s}{\tiny
Saddle cobordism}
 \epsfig{figure=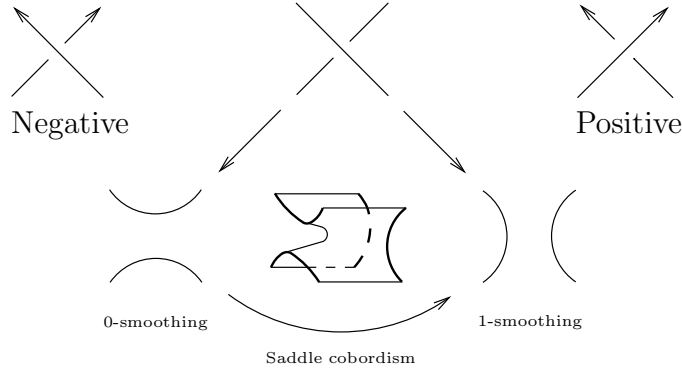}  \caption{Two smoothings of a crossing.} \label{smoothing.figure}
\end{center}
\end{figure}

\begin{itemize}
\item If $D$ has no crossings (therefore a flat tangle),
$\overline{\cF}(D)$ is just the complex
$$0\rightarrow \cF(D)\rightarrow 0,$$
where $\cF(D)$, sitting in cohomological degree zero, is the
bimodule associated to the flat tangle $D$.

\item If $D$ has one crossing, consider the complex
$\overline{\cF}(D)$ of $(A_m^{k,l},A_n^{k,l})$-bimodules
   \begin{equation*}
     0 \to \cF(D(0)) \stackrel{\mbox{\scriptsize{$\partial$}}}{\rightarrow}
\cF(D(1))\{-1\} \to
     0 \label{complex.equation}
   \end{equation*}
where $D(i), i=0,1$ denotes the $i$-smoothing of the crossing,
$\partial$ is induced by the ``saddle'' cobordism (see
figure~\ref{smoothing.figure}), and $\cF(D(0))$ sits in the
cohomological degree zero.

\item To a diagram with $t+1$ crossings we associate the total
complex $\overline{\cF}(D)$ of the bicomplex
   \begin{equation*}
   0 \to \cF(D(c_0)) \stackrel{\mbox{\scriptsize{$\partial$}}}{\rightarrow}
    \cF(D(c_1))\{-1\} \to 0 \label{complex1.equation}
   \end{equation*}
where $D(c_i), i=0,1$ denotes the $i$-smoothing of a crossing $c$
of $D.$

\item Finally, define $\cF(D)$ to be $\overline{\cF}(D)$ shifted
by $[x(D)]\{ 2x(D)-y(D)\}$, where $x(D)$ counts the number of
negative crossings and $y(D)$ counts the number of positive
crossings (see figure~\ref{smoothing.figure}).
\end{itemize}

\begin{figure}[ht!]
\begin{center}
\psfrag{pd}{\tiny $+$}\psfrag{nd}{\tiny
$-$}\psfrag{s1}{\scriptsize
$\{1\}$}\psfrag{s2}{\scriptsize$\{2\}$}\psfrag{s3}{\scriptsize$\{3\}$}
\psfrag{F}{$\cF$}\psfrag{d}{$\oplus$}\psfrag{e}{$:$}\psfrag{p}{$\partial$}
\psfrag{h}{\scriptsize Homological
grading:}\psfrag{-3}{\scriptsize$-3$}\psfrag{-2}{\scriptsize$-2$}\psfrag{-1}{\scriptsize$-1$}\psfrag{0}{\scriptsize$0$}
 \epsfig{figure=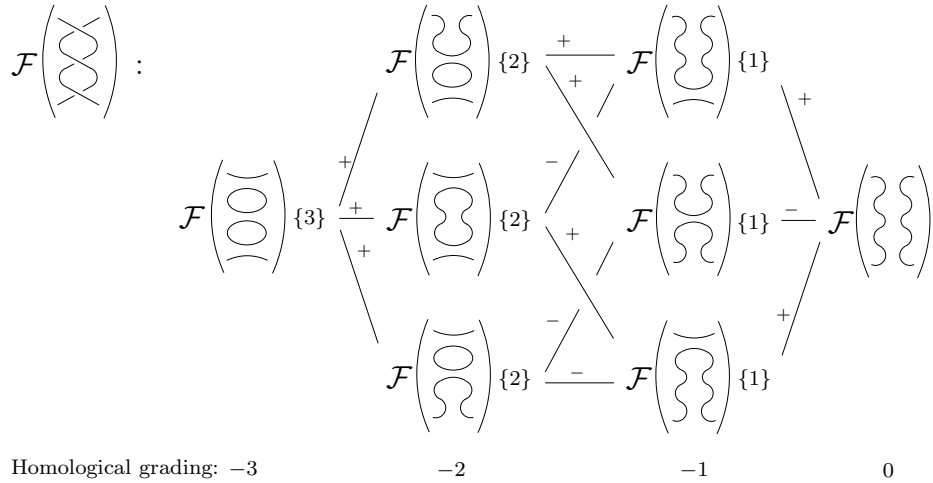}\caption{A total complex associated to a $(2,2)$-tangle.} \label{cube.figure}
\end{center}
\end{figure}

Figure~\ref{cube.figure} shows a complex of bimodules associated
to a $(2,2)$-tangle. Each arrow is induced by the saddle cobordism
and the sign on each arrow indicates the sign of each map in the
total complex.

\begin{theorem} If $D_1$ and $D_2$ are two diagrams of an oriented
$(m,n)$-tangle $L$, the complexes $\cF(D_1)$ and $\cF(D_2)$ of
graded $(A_{m}^{k,l},A_{n}^{k,l})$-bimodules are chain homotopy
equivalent.
\end{theorem}

The following proposition is a special case of the more general
theorem \ref{center.theorem} in section $3$.

\begin{prop} \label{DegreeZeroCentral.proposition}
The only invertible degree zero central elements in $A_n^{k,l}$
are $\pm 1$
\begin{equation*}
  Z^*_0(A_n^{k,l})\cong \{\pm 1\}.
\end{equation*}
\end{prop}

We now extend our invariant to oriented tangle cobordisms. Let $S$
be a movie presentation of a cobordism between two
$(m,n)$-tangles. $S$ is thus a sequence of Reidemeister moves and
handle moves. Each consequent pair of frames corresponds to a
homomorphism which is an isomorphism for each Reidemeister move,
and is induced by $\iota$, $\varepsilon$, $m$, or $\Delta$ for
each handle move. The composition of these homomorphisms gives us
a homomorphism $$\cF(S): \cF(D) \longrightarrow \cF(D'),$$ where
$D$ and $D'$ are the first and the last frame in the movie $S.$ It
follows from proposition \ref{DegreeZeroCentral.proposition} that
$\cF(S)=\pm \cF(S')$ if $S$ and $S'$ are two different
presentations of the same cobordism. To summarize, we have the
following:

\begin{prop} Complexes $\cF(T)$ of bimodules and homomorphisms
$\pm \cF(S)$ assigned to diagrams of tangle cobordisms assemble
into a projective $2$-functor from the $2$-category of oriented
tangle cobordisms to the $2$-category of natural transformations
between exact functors between homotopy categories of complexes of
graded $A_n^{k,l}$-modules.
\end{prop}

{\bf{Remark:}} The projective Grothendieck group
$K_p(A_n^{k,l}-\mbox{gmod})$ of the category of finitely-generated
graded projective $A_n^{k,l}$-modules is a free
$\mathbb{Z}[q,q^{-1}]$-module with a basis $[P_n^{k,l}(a)],$ $a\in
B_n^{k,l}.$ There is a natural way to identify
$K_p(A_n^{k,l}-\mbox{gmod})$ with a $\mathbb{Z}[q,q^{-1}]$-lattice
of $\mbox{Hom}(V_k\otimes V_l, V^{\otimes n})$
\begin{equation}\label{GrothendieckGroup.equation}
K_p(A_n^{k,l}-\mbox{gmod})\otimes_{\mathbb{Z}[q,q^{-1}]}\mathbb{C}\cong
\mbox{Hom}(V_k\otimes V_l, V^{\otimes n}),
\end{equation}
where $V$ is the fundamental $2$-dimensional representation of
$U_q(sl_2)$ and $V_i$ is the irreducible $(i+1)$-dimensional
representation of $U_q(sl_2)$. Under this isomorphism the basis
$[P_n^{k,l}(a)]$ gose to dual canonical basis of
$$\mbox{Hom}(V_k\otimes V_l, V^{\otimes n})\cong
\mbox{Inv}(V_k^{*}\otimes V_l^{*}\otimes V^{\otimes
n})\cong\mbox{Inv}(V_k\otimes V_l\otimes V^{\otimes n}).$$ Denote
by $\mathcal{K}(\mathcal{W})$ the category of bounded complexes of
objects of an abelian category $\mathcal{W}$ up to chain
homotopies. For each $(m,n)$-tangle $T$, it follows from
proposition~\ref{projectivity.proposition} that the complex of
bimodules $\cF(T)$ consists of right projective bimodules.
Therefore the tensor product with $\cF(T)$ is an exact functor
from $\mathcal{K}(A_n^{k,l}-\mbox{gmod})$ to
$\mathcal{K}(A_m^{k,l}-\mbox{gmod})$ which induces a homomorphism
$[\cF(T)]$ of $\mathbb{Z}[q,q^{-1}]$-modules
$$ K_p(A_n^{k,l}-\mbox{gmod}) \longrightarrow  K_p(A_m^{k,l}-\mbox{gmod}).$$
Direct computation shows that under the
isomorphism~(\ref{GrothendieckGroup.equation}) they give the
standard action of tangles on $\mbox{Inv}(V_k\otimes V_l\otimes
V^{\otimes n})$.

\section{The center of $A_n^{k,l}$}
We prove in this section that the center of the ring $A_n^{k,l}$
is isomorphic to the cohomology ring of
$\mathcal{B}_{\sigma_1,\sigma_2}$. Recall that
$\mathcal{B}_{\sigma_1,\sigma_2}$ denotes the Springer variety of
complete flags in $\mathbb{C}^n$ stabilized by a fixed nilpotent
operator with two Jordan blocks of sizes $\sigma_1$ and $\sigma_2$
respectively. Following from Khovanov's construction in
\cite{Khovanov:Crossingless}, we introduce the space
$\widetilde{S}$ and use it as a bridge to link the center of
$A_n^{k,l}$ and the cohomology rings of Springer varieties.
Without loss of generality, we assume throughout this section that
$n\geq m$, $n+m\equiv 0$ mod $2$, and $0\leq l-k \leq n$ (note
that $A_n^{k,l}$ is trivial if $l-k>n$ ). The proofs in this
section rely heavily on \cite{Khovanov:Crossingless}.

\begin{theorem}\label{center.theorem}
The center of $A_n^{k,l}$ is isomorphic to the cohomology ring of
$\mathcal{B}_{\sigma_1,\sigma_2}$
$$Z(A_n^{k,l})\cong H^*(\mathcal{B}_{\sigma_1,\sigma_2}),$$
where $\sigma_1=\frac{n+l-k}{2}$ and $\sigma_2=\frac{n-l+k}{2}$.
\end{theorem}

Denote by $S$ the $2$-sphere $S^2$ and let $p$ be the north pole
of $S$. Let $S^{\times n}$ be the direct product of $n$ spheres
$$S^{\times n}\define \underbrace{S\times S\times \cdots \times S}_{n}. $$
Label the $n$ free points of $B_n^{k,l}$ by $1,2,\cdots,n$ from
left to right. For each $a\in B_n^{k,l}$ define a submanifold
$S_a\in S^{\times n}$ consisting of sequences $(x_1,\cdots x_n)$,
$x_i\in S$, such that $x_i=x_j$ whenever $(i,j)$ is a type I arc
in $a$, and $x_s=p$ if $s$ is connected to a platform. Let
$\widetilde{S}_n^{k,l}$ be the subspace of $S^{\times n}$ which is
the union of all $S_a$
$$\widetilde{S}_n^{k,l}\define \bigcup_{a\in B_n^{k,l}} S_a.$$
When there is no confusion we write $\widetilde{S}$ instead of
$\widetilde{S}_n^{k,l}$.

Note that the cohomology ring of $S_a$ is isomorphic to the ring
$_a(A_n^{k,l})_a$, and the cohomology ring of $S_a \cap S_b$,
viewed as abelian group, is isomorphic to $_a(A_n^{k,l})_b$. These
observations lead to the following:

\begin{theorem}\label{centerS.theorem}
The center of $A_n^{k,l}$ is isomorphic to the cohomology ring of
$\widetilde{S}$
$$Z(A_n^{k,l})\cong H^*(\widetilde{S},\mathbb{Z}).$$
\end{theorem}
\begin{proof} Denote by $H(Y)$ the cohomology ring of the space
$Y$ with integer coefficients. As noted above, we have
$H(S_a)\cong {_a(A_n^{k,l})_a}$ and $H(S_a \cap S_b)\cong
{_a(A_n^{k,l})_b}$. The second isomorphism allows us to make
$_a(A_n^{k,l})_b$ into a ring with unit $1\define 1^{s} \in
\cF(W(b)a)\cong \cA^{\otimes s}$.

We have natural ring homomorphisms induced by inclusions
$$\psi_{a;a,b}:H(S_a)\rightarrow H(S_a \cap S_b), \hspace{0.2 in} \psi_{b;a,b}:H(S_b)\rightarrow H(S_a \cap S_b),$$
and also
$$\gamma_{a;a,b}:{_a(A_n^{k,l})_a}\rightarrow {_a(A_n^{k,l})_b}, \hspace{0.2 in} \gamma_{b;a,b}:{_b(A_n^{k,l})_b}\rightarrow {_a(A_n^{k,l})_b},$$
which are given by $x\mapsto x {_a 1_b}$ and $x\mapsto {_a 1_b}x$.
Assemble all these together we get a commutative diagram of ring
homomorphisms
 \begin{equation*}
    \begin{CD}
      H(\widetilde{S}) @>\tau>>  \mathrm{Eq}(\psi)    @>>>  {\mathop{\prod}\limits_{a}}H(S_a) @>{\psi}>>
    {\mathop{\prod}\limits_{a\not= b}} H(S_a\cap S_b) \\
      &&    @VV{\cong}V     @VV{\cong}V                @VV{\cong}V     \\
     Z(A_n^{k,l})@>\cong>> \mathrm{Eq}(\gamma) @>>> {\mathop{\prod}\limits_{a}}
     {_a(A_n^{k,l})_a} @>{\gamma}>>
    {\mathop{\prod}\limits_{a\not= b}}  {_a(A_n^{k,l})_b}
    \end{CD}
   \end{equation*}
where
\begin{equation*}
\psi = \sum_{a\not= b}(\psi_{a;a,b}+ \psi_{b;a,b})\hspace{0.1in}
 \mathrm{and} \hspace{0.1in}
 \gamma = \sum_{a\not= b}(\gamma_{a;a,b}+ \gamma_{b;a,b}).
\end{equation*}
$\mathrm{Eq}(\alpha)$ is the \emph{equalizer} of the map $\alpha$
(see \cite{Khovanov:Crossingless}). For example,
$\mathrm{Eq}(\psi)$ is a subring of ${\mathop{\prod}\limits_{a}}
H(S_a)$ consisting of $\times_a h_a$ such that if $h_a\in
H(S_a)$ and $h_b\in H(S_b)$ then their images in $H(S_a\cap S_b)$
under $\psi$ are equal.

For $\forall x\in A_n^{k,l}$ write $x$ as $\sum_{a,b\in B_n^{k,l}}
{_a x_b}$. Assuming $x$ is central we have $x 1_b 1_a=1_b x 1_a =
{_a x_b}$. Therefore ${_a x_b}=0$ if $a\neq b$. So $x=\sum_{a} {_a
x_a}$ is central if and only if $({_a x_a})({_a 1_b})=({_a
1_b})({_b x_b})$, which means $Z(A_n^{k,l})\cong
\mathrm{Eq}(\gamma)$. The ring homomorphism
$H(\widetilde{S})\rightarrow {\mathop{\prod}\limits_{a}} H(S_a)$
factors through $\mathrm{Eq}(\psi)$. To prove theorem 2 it
suffices to show that $\tau$ is an
isomorphism.\\

For $a,b\in B_n^{k,l}$ write $a\rightarrow b$ if there exists a
``horizontal merging'' of two arcs (see
figure~\ref{merge_hv.figure}).
\begin{figure}[ht!]
\begin{center}
\psfrag{h}{\small Horizontal merging}\psfrag{v}{\small Vertical
merging}\epsfig{figure=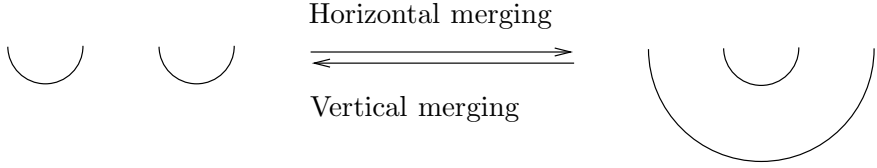} \caption{Horizontal and
vertical mergings of two arcs.} \label{merge_hv.figure}
\end{center}
\end{figure}
Introduce a partial order on $B_n^{k,l}$ by setting $a\prec b$ if
there is a chain of arrows $a\rightarrow a_1\rightarrow \cdots
\rightarrow a_m\rightarrow b$. Extend this partial order
arbitrarily to a total order $<$ on $B_n^{k,l}$. See
figure~\ref{ArrowRelations.figure} for arrow relations and
ordering of $B_5^{0,1}$ ($a_i<a_j$ iff $i<j$).

\begin{figure}[ht!]
\begin{center}
\psfrag{a1}{\tiny $a_1$}\psfrag{a2}{\tiny $a_2$}\psfrag{a3}{\tiny
$a_3$}\psfrag{a4}{\tiny $a_4$}\psfrag{a5}{\tiny $a_5$}
\epsfig{figure=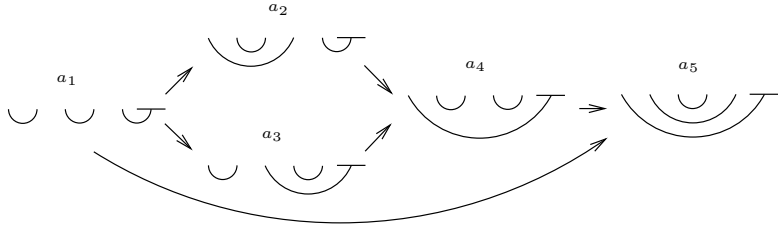} \caption{Arrow relations and ordering of
$B_5^{0,1}$.} \label{ArrowRelations.figure}
\end{center}
\end{figure}

We would like to construct a cell decomposition of $S_a$.
Associate a decorated graph $\Gamma$ to $a\in B_n^{0,m}$ as
follows (see figure~\ref{graphofa.figure} for an example):
\begin{itemize}
\item Each type I arc $x_i$ in $a$ corresponds to a ``hollow''
vertex $i$ in $\widetilde{\Gamma}$.

\item Each type II arc $x_j$ in $a$ corresponds to a ``solid''
vertex $j$ in $\widetilde{\Gamma}$.

\item Two vertices $i,j$ in $\widetilde{\Gamma}$ are connected by
an edge iff the result of merging $x_i$ and $x_j$ ``vertically''
still lies in $B_n^{0,m}$.

\item $\Gamma$ is obtained from $\widetilde{\Gamma}$ by
contracting all edges with two ``solid'' ends.

\item Mark a vertex in each connected component of $\Gamma$
without ``solid'' vertices.

\end{itemize}

\begin{figure}[ht!]
\begin{center}
\psfrag{m}{\tiny $m$} \psfrag{a}{$a$}
\psfrag{G1}{$\widetilde{\Gamma}$} \psfrag{G2}{$\Gamma$}
\epsfig{figure=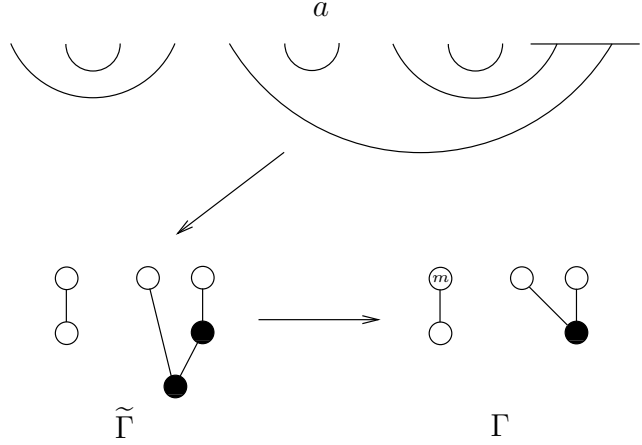} \caption{An element $a\in
\widetilde{B}_{10}^{0,2}$ and its associated graph $\Gamma$.}
\label{graphofa.figure}
\end{center}
\end{figure}

Let $E$ be the set of edges, $M$ be the set of marked points, and
$I$ be the set of vertices in $\Gamma$. For each $J\subset
(E\sqcup M)$ let $c(J)$ be the subset of $S^{\times I}$ consisting
of points $\{ y_i\}_{i\in I}, y_i\in S$ such that
 \[
 \begin{array}{ll}
  y_i = y_j       & \mbox{if $(i,j)\in J$}, \\
  y_i \not= y_j   & \mbox{if $(i,j)\notin J$}, \\
  y_i = p         & \mbox{if $i\in M\cap J$}, \\
  y_i = p         & \mbox{if $i$ is ``solid''}, \\
  y_i \not= p     & \mbox{if $i\in M, i\notin J$}.
 \end{array}
 \]
Clearly, $S^{\times I}= \sqcup_{J} c(J)$ and $c(J)$ is an open
cell of dimension $2(|I|-|J|)$. We thus obtain a decomposition of
$S_a$ into even dimensional cells.

\begin{lemma}\label{intersection.lemma}
$S_{<a}\cap S_a = (\cup_{b\rightarrow a} S_b)\cap S_a$, where $S_{< a}=\bigcup_{b< a}S_b$.
\end{lemma}

It follows from lemma \ref{intersection.lemma} and the above
construction that

\begin{lemma} \label{partition.lemma}
The cell decomposition constructed above restricts to a cell
decomposition of $S_a\backslash S_{<a}$, which is a union of cells
$c(J)$ such that $J\cap E= \emptyset$.
\end{lemma}

We thus obtain a \emph{cell partition} of $\widetilde{S}$ by
adding cells in $S_a \backslash S_{<a}$ following the total order.
Since there are only even-dimensional cells in the partition, the
rank of $H(\widetilde{S})$ is equal to the number of cells.

By induction on $a$ with respect to the total order $<$ we
get the following (see \cite{Khovanov:Crossingless}):

\begin{prop} \label{prop-exact}
$S_{\leq a}$ has cohomology in even degrees only and the following
sequence is exact
\begin{equation}\label{equation-exact}
 0 \rightarrow H(S_{\le a}) \stackrel{\phi}{\rightarrow} \oplusop{b\le a} H(S_b)
  \stackrel{\psi^-}{\rightarrow}
 \oplusop{b< c\le a} H(S_b\cap S_c),
\end{equation}
where $\phi$ is induced by inclusions $S_b\subset S_{\le a},$
while
\begin{equation*}
\psi^- \define \sum_{b<c\le a} (\psi_{b,c}-\psi_{c,b}),
\end{equation*}
where
\begin{equation*}
\psi_{b,c}: H(S_b) \rightarrow H(S_b\cap S_c)
\end{equation*}
is induced by the inclusion $(S_b\cap S_c )\subset S_b.$
\end{prop}

When $a$ is maximal with respect to the total order $<$, the exact
sequence (\ref{equation-exact}) becomes
\begin{equation*}
 0 \rightarrow H(\widetilde{S}) \stackrel{\phi}{\rightarrow} \oplusop{b} H(S_b)
  \stackrel{\psi^-}{\rightarrow}
 \oplusop{b< c} H(S_b\cap S_c),
\end{equation*}
which means that $H(\widetilde{S})$ is isomorphic to the equalizer
of $\psi$.
\end{proof}

\begin{lemma}\label{oneplatform.lemma}
The center of $A_n^{k,l}$ is isomorphic to the center of
$A_n^{0,l-k}$
$$Z(A_n^{k,l})\cong Z(A_n^{0,l-k}).$$
\end{lemma}
\begin{proof} It follows from the definition of $\widetilde{S}$ that
$$\widetilde{S}\cong \bigcup_{a\in \widetilde{B}_n^{k,l}} S_a,$$
where $\widetilde{B}_n^{k,l}\subset B_n^{k,l}$ is the set of
crossingless matchings with all points on the left platform
connected to the right platform.
\begin{figure}[h!]
\begin{center}
\psfrag{bt}{$\widetilde{B}_3^{2,3}$:}
\psfrag{b}{$B_3^{0,1}$:}\epsfig{figure=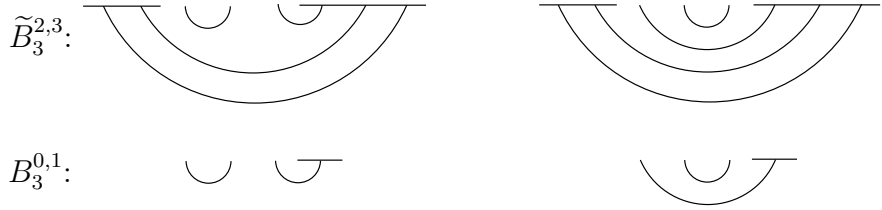}
\caption{Removing bottommost type II arcs.} \label{bt312.figure}
\end{center}
\end{figure}
On the other hand, since the bottommost type II arcs contribute
nothing, we have (see figure~\ref{bt312.figure})
$$\bigcup_{a\in\widetilde{B}_n^{k,l}}S_a\cong \bigcup_{a\in B_n^{0,l-k}}S_a.$$
Theorem $2$ and the above observations imply that
$$Z(A_n^{k,l})\cong H(\widetilde{S}_n^{k,l})\cong H(\widetilde{S}_n^{0,l-k})\cong Z(A_n^{0,l-k}).$$
\end{proof}

\begin{prop}\label{dimension.lemma}
$\widetilde{S}_n^{0,m}$ has {\scriptsize
$\sbinom{n}{\frac{n-m}{2}}$} cells in the partition constructed
above.
\end{prop}
\begin{proof} The proposition is equivalent to the statement that
$\widetilde{S}_{2s-k}^{0,k}$ has {\scriptsize
$\sbinom{2s-k}{s-k}$} cells. Fix the total number of points $2s$
(including marked and free). Induction on the size of the
right platform $k$.

 Induction base $k=0$ is proved in \cite{Khovanov:Crossingless}, lemma
 4.1. Assuming the statement is true up to $k$, it suffices to prove
 that extending the size of the platform by $1$ eliminates {\scriptsize
$\sbinom{2s-k}{s-k}$} - {\scriptsize $\sbinom{2s-k-1}{s-k-1}$}
cells in $\widetilde{S}_{2s-k}^{0,k}$. If we label the $2s$ points
by $1,2,3,\cdots, 2s$ from right to left, the vanished cells are
exactly those in $S_a$ where $a$ has a type II arc connecting $k$
and $k+1$ (see figure~\ref{vanishedcells.figure}). Denote the set
of those $a$ by $a(k,k+1)$.
\begin{figure}[ht!]
\begin{center}
\psfrag{1}{\tiny $1$}\psfrag{k-1}{\tiny $k-1$}\psfrag{k}{\tiny
$k$}\psfrag{k+1}{\tiny $k+1$}\psfrag{ik-1}{\tiny
$i_{k-1}$}\psfrag{i1}{\tiny $i_1$}\psfrag{i1+1}{\tiny
$i_1+1$}\psfrag{2s}{\tiny $2s$}\psfrag{cd}{\tiny
$\cdots$}\psfrag{B}{\tiny
$B_{2s-i_1}^{0,0}$}\psfrag{in}{\begin{rotate}{-90}$\in$\end{rotate}}
\epsfig{figure=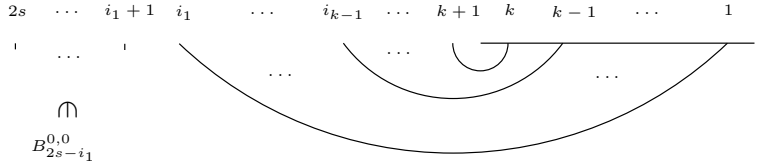} \caption{A generic element in
$a(k,k+1)$.} \label{vanishedcells.figure}
\end{center}
\end{figure}
We thus get the following formula:
\begin{equation*}
\vartheta(\bigcup_{a(k,k+1)\in B_{2s-k}^{0,k}} S_{a(k,k+1)})=
\sum_{i=0}^{s-k}(i+1)C_i [(\frac{1-\sqrt{1-4x}}{2x})^{k-1}]_{i},
\end{equation*}
where $\vartheta(X)$ denotes the number of cells in $X$, $C_i$
denotes the $n$-th Catalan number, and $[f(x)]_i$ is the
coefficient of $x^i$ in the Maclaurin expansion of $f(x)$. Recall
that $\frac{1-\sqrt{1-4x}}{2x}$ is the generating function of
$\{C_i\}$. The factor $(i+1)C_i$ corresponds to the number of
cells outside the bottommost type II arc and
$[(\frac{1-\sqrt{1-4x}}{2x})^{k-1}]_{i}$ is equal to the total
number of crossingless matchings inside the bottommost type II
arc. Note that the arcs inside the bottommost type II arc
contribute nothing to the cell structure, see
figure~\ref{graphofa.figure}. After simplifying the above formula
we get
\begin{equation*}
\vartheta(\bigcup_{a(k,k+1)\in B_{2s-k}^{0,k}}
S_{a(k,k+1)})=[\frac{1}{\sqrt{1-4x}}(\frac{1-\sqrt{1-4x}}{2x})^{k-1}]_{s-k}.
\end{equation*}
By induction on $s$ and $k$ it's easy to prove that
\begin{equation*}
[\frac{1}{\sqrt{1-4x}}(\frac{1-\sqrt{1-4x}}{2x})^{k-1}]_{s-k} =
{\scriptsize \sbinom{2s-k}{s-k}} - {\scriptsize
\sbinom{2s-k-1}{s-k-1}}
\end{equation*}
and the proposition follows.
\end{proof}

\begin{prop}\label{prop.ConciniProcesi} \cite{ConciniProcesi}
The cohomology ring of $\mathcal{B}_{\frac{n+m}{2},\frac{n-m}{2}}$
has dimension ${\scriptsize \sbinom{n}{\frac{n-m}{2}}}$ and is
  isomorphic to the quotient ring of $R=\mathbb{Z}[X_1,\dots, X_{n}]$ by the ideal
  $R_1$ generated by $e_k(I)$ for all $k+|I|=n+1$,  $X_I$ for all
  $|I|=\frac{n-m}{2}+1$, and $X_i^2$ for $i\in [1,n]$, where
  \begin{equation*}
  I\subset \{1,2,\cdots,n\},\ \ \ \  X_I=\prod_{i\in I}X_i,\ \ \ \ e_k(I) = \sum_{|J|=k, J\subset I} X_J.
 \end{equation*}
\end{prop}

We now prove the main theorem of this section.

\emph{Proof of theorem~\ref{center.theorem}.} It follows from
theorem~\ref{centerS.theorem} and lemma~\ref{oneplatform.lemma}
that to prove theorem~\ref{center.theorem} it suffices to show
that $H(\widetilde{S}_n^{0,m})\cong
H(\mathcal{B}_{\frac{n+m}{2},\frac{n-m}{2}})$. Denote by $X$ a
generator of $H^2(S)$. We have the following maps
\begin{equation*}
    \begin{CD}
      \widetilde{S}_n^{0,m} @>\iota>>  S^{\times n}    @>{\psi_i}>>
      S',
    \end{CD}
\end{equation*}
where $\iota$ is the inclusion and $\psi_i$ is the projection onto the $i$-th component.
Define $X_i\in \widetilde{S}_n^{0,m}$ to be the pull back of $X$ under the map
$\psi_i\circ\iota$
\begin{equation*}
X_i\define (-1)^i \iota^* \circ \psi_i^* (X).
\end{equation*}
It's obvious that those $\{X_i\}$ generate $H(\widetilde{S}_n^{0,m})$. It
follows from proposition~\ref{prop.ConciniProcesi} and
proposition~\ref{dimension.lemma} that to prove the theorem we only need
to verify the following relations:

 \begin{eqnarray}
  X_i^2 & = & 0, \hspace{0.2in} i\in [1,n]; \\
  X_I   & = & 0, \hspace{0.2in} |I|=\frac{n-m}{2}+1;\\
  e_k(I)   & = & 0, \hspace{0.2in} k+|I|=n+1. \label{relation.equation}
 \end{eqnarray}

The first two relations are obvious. Consider the map $i^*_a:
H(\widetilde{S}_n^{0,m})\rightarrow H(S_a)$ induced by the
inclusion $i_a: S_a \hookrightarrow \widetilde{S}_n^{0,m}$. Since
$\sum_{a} i_a^*(H(\widetilde{S}_n^{0,m}))\rightarrow
\oplus_{a}H(S_a)$
 is an inclusion, (\ref{relation.equation}) will follow once we verify that that
\begin{equation}\label{relation1.equation}
\sum_{|J|=k, J\subset I} i_a^*(X_J)=0,\hspace{0.2in} k+|I|=n+1
\end{equation}
for all $a\in B_n^{0,m}$.

Fix any $a\in B_n^{0,m}$ and $I\subset \{1,2,\cdots,n\}$. Since
$n-|I|=k-1$ there exists at most $k-1$ type I arcs where $I$
intersects with each one at only one point. Therefore, for each
$J\subset I$ such that $|J|=k$,
 $J$ must either contain an end point of a
type II arc or contain an end point of a type I arc $(p_1,p_2)$
such that $\{p_1,p_2\}\in I$. If $J$ contains an end point of a
type II arc then $i_a^*(X_J)=0$. For a type I arc $(p_1,p_2)$,
because of the term $(-1)^i$ in the definition of $X_i$ and the
fact that $p_1+p_2$ is odd, we have $i_a^* (X_{p_1} X_{p_2})=0$
and $i_a^*(X_{p_1} + X_{p_2})=0$. Therefore
\begin{equation*}
  \sum_{J\subset I, |J|=k, \{p_1,p_2\}\cap J \neq \emptyset} i_a^* (X_J)
  =0.
\end{equation*}
For the remaining terms in the summation in
(\ref{relation1.equation}), pick another type I arc and repeat the
above process. After finitely many reductions we can get the
relation (\ref{relation1.equation}).\hfill$\square$

\section{Categorification of level two representations of quantum $sl_N$}
\subsection{Level two representations of quantum $sl_N$}
Let $W$, $\wedge^2 W,\cdots$, $\wedge^{N-1}W$ be the irreducible
representations of $U_q(sl_N)$ with highest weights $\omega_1,\
\omega_2,\ \cdots,\ \omega_{N-1}$ respectively, where
$\omega_i=L_1+\cdots+L_{i}$ and the $L_j$'s are the fundamental
weights. A level two representation $V$ of $U_q(sl_N)$ is an
irreducible representation with the highest weight
$\lambda=\omega_s+\omega_{s+k}$ for some $k\geq 0$. Fix $V$ for
the rest of this paper. $V$ decomposes into weight spaces $V=
\oplusop{\mu} V_{\mu}$. Following \cite{HK:LevelTwo}, we call
$\mu$ admissible if $\mu$ appears in the weight space
decomposition of $V$. A weight $\mu$ is admissible if and only if
it can be written as the sum $\mu_1 L_1+\mu_2 L_2+\cdots +\mu_N
L_N$ such that
\begin{itemize}
\item $0\leq \mu_i \leq 2,$ for $\forall 1\leq i\leq N$

\item $\sum_{i=1}^{N}\mu_i=2s+k$

\item $\mu_1+\cdots+\mu_i\leq \lambda_1+\cdots+\lambda_i$, for
$\forall 1\leq i\leq N$,
\end{itemize}
where $\lambda_i$ is the coefficient of $L_i$ in the decomposition
$$\lambda=\omega_s+\omega_{s+k}=(L_1+\cdots+\L_s)+(L_1+\cdots+\L_{s+k}).$$
For each admissible weight $\mu$ let $m(\mu)$ be the number
of $1$'s in the sequence $(\mu_1,\dots,\mu_N)$. The dimension of
$V_{\mu}$ is then determined by $m(\mu)$.

Recall that $U_q(sl_N)$ is defined to be the algebra generated by
$E_i$, $F_i$, $K_i$, and $K_i^{-1}$ for $1\leq i \leq N-1$ with
relations
\begin{equation} \label{q-rel}
\begin{array}{l}
 K_i K_i^{-1} = 1 = K_i^{-1} K_i, \\
 K_i K_j = K_j K_i, \\
 K_i E_j = q^{c_{i,j}} E_j K_i, \\
 K_i F_j = q^{-c_{i,j}} F_j K_i, \\
 E_i F_j - F_j E_i = \delta_{i,j} \frac{K_i - K_i^{-1}}{q-q^{-1}}, \\
 E_i E_j = E_j E_i \ \ \mathrm{if} \ \ |i-j|>1, \\
 F_i F_j = F_j F_i \ \ \mathrm{if} \ \ |i-j|>1, \\
 E_i^2 E_{i\pm 1} - (q+q^{-1}) E_i E_{i\pm 1}E_i + E_{i\pm 1}E_i^2 = 0, \\
 F_i^2 F_{i\pm 1} - (q+q^{-1}) F_i F_{i\pm 1}F_i + F_{i\pm 1}F_i^2 = 0.
 \end{array}
\end{equation}
$E_i$ acts on $V$ by sending weight space $V_\mu$ to
$V_{\mu+\epsilon_i}$ and $F_i$ maps $V_\mu$ to
$V_{\mu-\epsilon_i}$, where
$\epsilon_i=(\underbrace{0,\cdots,0}_{i-1},1,-1,0,\cdots,0)$.

\subsection{Semi-standard tableaux and arc rings}
\begin{figure}[ht!]
\begin{center}
\psfrag{lambda}{\small $\lambda=(2,2,1,0)$}\psfrag{mu}{\small
$\mu=(1,1,2,1)$}\psfrag{1}{\small $1$}\psfrag{2}{\small
$2$}\psfrag{3}{\small $3$}\psfrag{4}{\small $4$}
 \epsfig{figure=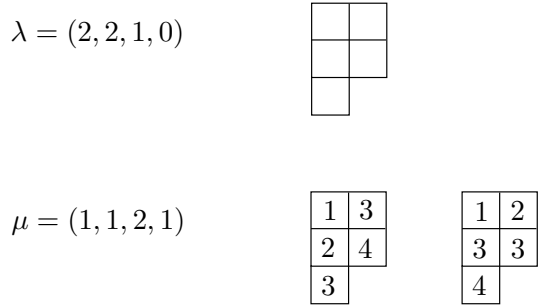}\caption{Semi-standard tableaux.} \label{YoungTableau.figure}
\end{center}
\end{figure}

We give in this section an explicit bijection between
semi-standard tableaux and crossingless matchings with one
platform. First recall the definition of semi-standard tableaux.
For any $\lambda=(\lambda_1,\cdots,\lambda_{N-1},0)$ in the weight
lattice of $U_q(sl_N)$, there exists an irreducible representation
$V_{\lambda}$ with the highest weight $\lambda$. Weight
$\mu=(\mu_1,\cdots,\mu_N)$ appears in the weight space
decomposition of $V$ if and only if $\mu$ is admissible. The
dimension of the weight space $V_{\lambda}(\mu)$ equals to the
number of ways one can fill the Young diagram corresponding to
$\lambda$ with $\mu_1\ 1's$, $\mu_2\ 2's$, $\cdots$, and $\mu_N\
N's$ in such a way that each column is strictly increasing and
each row is non-decreasing. Each such filling is called a
semi-standard tableau (see figure~\ref{YoungTableau.figure}).

\begin{figure}[ht!]
\begin{center}
\psfrag{1}{\small 1}\psfrag{2}{\small 2}\psfrag{3}{\small
3}\psfrag{4}{\small 4}\psfrag{5}{\small 5}\psfrag{6}{\small
6}\psfrag{7}{\small 7}\psfrag{8}{\small 8}\psfrag{9}{\small
9}\psfrag{cd}{\small $\cdots$}\psfrag{l}{\small $M_{\mu}\bigcap
T^i_{\mu}(l)=\{1,2,3,6,7\}$}\psfrag{r}{\small $M_{\mu}\bigcap
T^i_{\mu}(r)=\{4,5,8\}$}
 \epsfig{figure=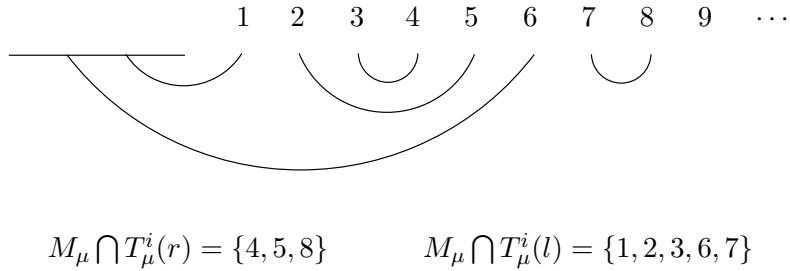}\caption{An element in $B_8^{2,0}(M_{\mu})$.} \label{MatchingExample.figure}
\end{center}
\end{figure}

The Young diagram $Y_{\lambda}$ corresponding to the level two
representation $V$ with the highest weight
$\lambda=\omega_s+\omega_{s+k}$ has two columns of length $s+k$
and $s$ respectively. Fix an admissible weight $\mu$. Let
$M_{\mu}$ be the set of $i$'s such that $\mu_{i}=1$
$$M_{\mu}\define \{1\leq i\leq N|\mu_i=1\},$$
and $N_{\mu}$ be the set of $i$'s such that $\mu_{i}=2$
$$N_{\mu}\define \{1\leq i\leq N|\mu_i=2\}.$$
Note that $|M_{\mu}|=m(\mu)$. Let $T_{\mu}$ be the set of
semi-standard tableaux of $Y_{\lambda}$ corresponding to $\mu$.
For each semi-standard tableau $T^i_{\mu}\in T_{\mu}$, let
$T^i_{\mu}(r)$ and $T^i_{\mu}(l)$ be the set of numbers on the
right and left column of $T^i_{\mu}$ respectivily. Write $M_{\mu}$
as an ordered sequence $\{a_1,a_2,\cdots,a_{m(\mu)}\}$. Assume
that $\{a_{i_1},a_{i_2},\cdots,a_{i_t}\}=M_{\mu}\bigcap
T^i_{\mu}(r)$. Consider all integer points $\{1,2,3,\cdots\}$
lying on the x-axis. Put a platform on the x-axis to the left of
all points (see figure~\ref{MatchingExample.figure}). First draw
an arc in the lower half plane connecting $a_{i_1}$ with the first
point in $M_{\mu}$ to its left which is not connected to any
point. That point always exists and lies in $T^i_{\mu}(l)$ since
$T^i_{\mu}$ is semi-standard. Repeat the above step for $a_2,a_3,
\cdots$ in order until each point in $M_{\mu}\bigcap T^i_{\mu}(r)$
is connected to some point. Finally, connect the remaining free
points in $M_{\mu}\bigcap T^i_{\mu}(l)$ to the platform by arcs in
the unique way that no two arcs intersect.
The resulting graph is a crossingless matching among the points in $M_{\mu}$
with one platform. Denote by $B_{m(\mu)}^{k,0}(M_{\mu})$ the set of all such elements.
Note that $B_{m(\mu)}^{k,0}(M_{\mu})\cong B_{m(\mu)}^{k,0}$. The map from
$T_{\mu}$ to $B_{m(\mu)}^{k,0}(M_{\mu})$ is denoted by $\varphi_{\mu}$.

\begin{figure}
\begin{center}
\psfrag{lambda}{\small $\lambda=(2,2,1,0,0)$}\psfrag{yd}{\small
Young diagram}\psfrag{mu1}{\small
$\mu=(1,1,1,1,1)$}\psfrag{mu2}{\small
$\mu=(2,1,1,1,0)$}\psfrag{1}{\small $1$}\psfrag{2}{\small
$2$}\psfrag{3}{\small $3$}\psfrag{4}{\small $4$}\psfrag{5}{\small
$5$}\psfrag{phi}{\small $\ \ \varphi$}\psfrag{psi}{\small $\ \
\psi$}\psfrag{relabel}{\small relabel} \psfrag{p1}{\small
$1$}\psfrag{p2}{\small $2$}\psfrag{p3}{\small
$3$}\psfrag{p4}{\small $4$}\psfrag{p5}{\small
$5$}\psfrag{mmu1}{\small
$M_{\mu}=\{1,2,3,4,5\}$}\psfrag{mmu2}{\small $M_{\mu}=\{2,3,4\}$}
 \epsfig{figure=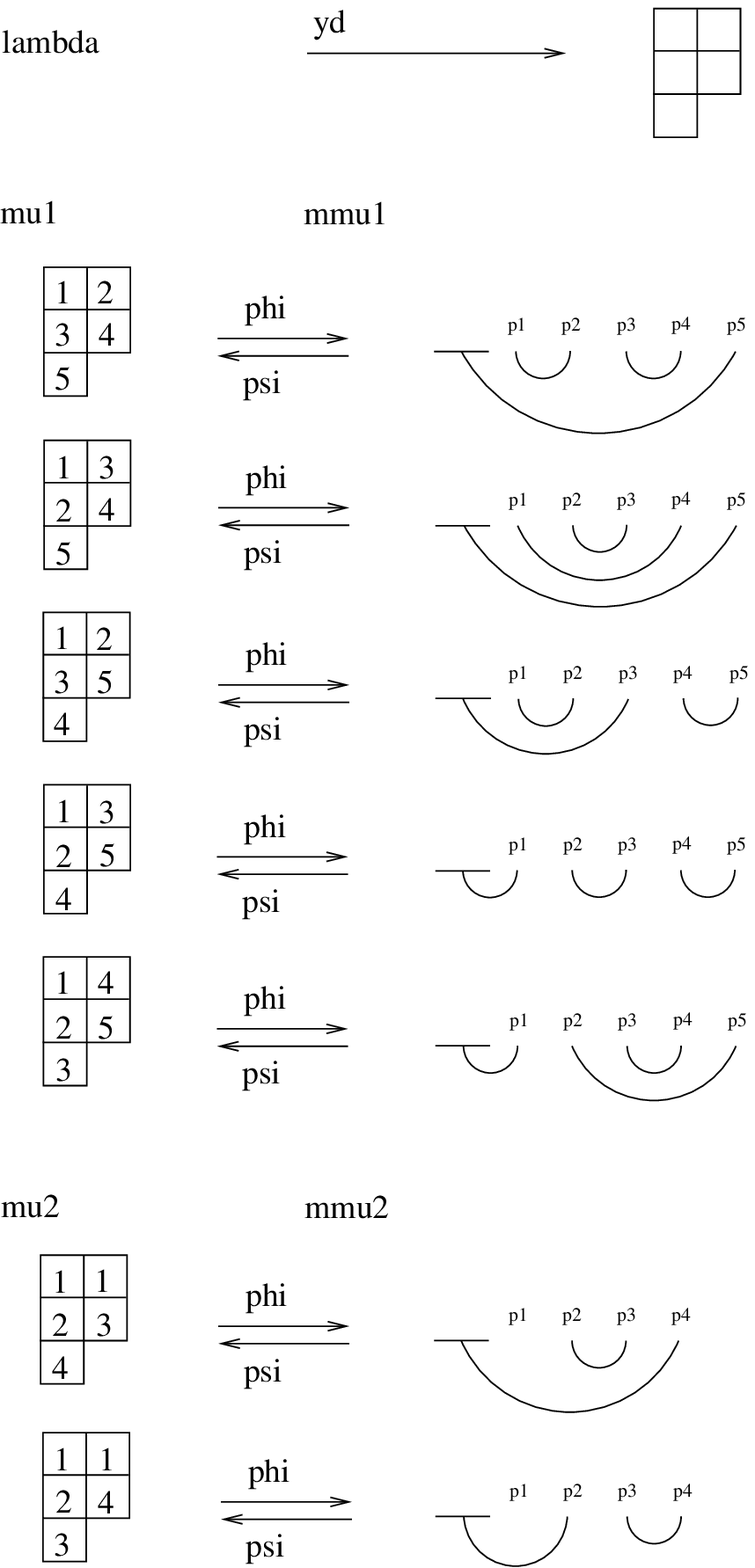}
 \caption{Bijection between semi-standard tableaux and crossingless matchings with one platform.}
 \label{YoungTableauBijection.figure}
\end{center}
\end{figure}

Conversely, for any $a\in B_{m(\mu)}^{k,0}(M_{\mu})$ with $t$ type
I arcs $c_1,c_2,\cdots,c_t$, let $R_a$ be the set of right end
points of all $c_i$. A semi-standard tableau of $Y_{\lambda}$ is
constructed by putting $R_a \bigcup N_{\mu}$ into the right column
and $(M_{\mu} \backslash R_a)\bigcup N_{\mu}$ into the left column
(both in increasing order). The map from $B_{m(\mu)}^{k,0}$ to
$T_{\mu}$ is denoted by $\psi_{\mu}$. It's easy to verify that
$\psi_{\mu}$ is indeed the inverse of $\varphi_{\mu}$. Thus we
have a bijection between $T_{\mu}$ and $B_{m(\mu)}^{k,0}$
\begin{equation}\label{bijection.equation}
T_{\mu} \leftrightmaps{50}{\mbox{\scriptsize{$\psi_{\mu}$}}}
{\mbox{\scriptsize{$\varphi_{\mu}$}}} B_{m(\mu)}^{k,0}.
\end{equation}
See figure~\ref{YoungTableauBijection.figure} for an example of
this bijection where $V$ is the level two representation of
$U_q(sl_5)$.

\subsection{Category $\mathcal{C}$ and exact functors}

Starting with $B_{m(\mu)}^{k,0}(M_{\mu})$, we repeat the definition
of $A_n^{k,l}$ and get $A_{m(\mu)}^{k,0}(M_{\mu})$, or simply
$A_{\mu}$. Note that $A_{\mu}\cong \mathbb{Z}$ when $m(\mu)=0$.
For an admissible weight $\mu$ define $\mathcal{C}_{\mu}$ to be
the category of finitely generated graded left $A_{\mu}$-modules.
By taking direct sum over all admissible $\mu$ we collect those
$\mathcal{C}_{\mu}$ into a single category $\mathcal{C}$
$$\mathcal{C}\define \oplusop{\mu} \mathcal{C}_{\mu}.$$
Note that when $k=0$ our category $\mathcal{C}_{\mu}$ is the same
as $\mathcal{C}(\mu)$ in \cite{HK:LevelTwo}.

\begin{figure}[ht!]
\begin{center}
\psfrag{1}{\small $1$}\psfrag{2}{\small $2$}\psfrag{0}{\small $0$}
\psfrag{mu1}{\small $\mu_1=(1,1,0,1,2,1)$}\psfrag{mu2}{\small
$\mu_1+\epsilon_3=(1,1,1,0,2,1)$}\psfrag{mu3}{\small
$\mu_2=(1,2,1,1,1,1)$}\psfrag{mu4}{\small
$\mu_2+\epsilon_3=(1,2,2,0,1,1)$}\psfrag{e3}{\small
$\varepsilon_3=(0,0,1,-1,0,0)$}\psfrag{E3}{\small
$\mathcal{E}_3^{\mu_1}=\cF$}\psfrag{E4}{\small
$\mathcal{E}_3^{\mu_2}=\cF$}\psfrag{p}{\LARGE
$)$}\psfrag{l}{\LARGE $($}
 \epsfig{figure=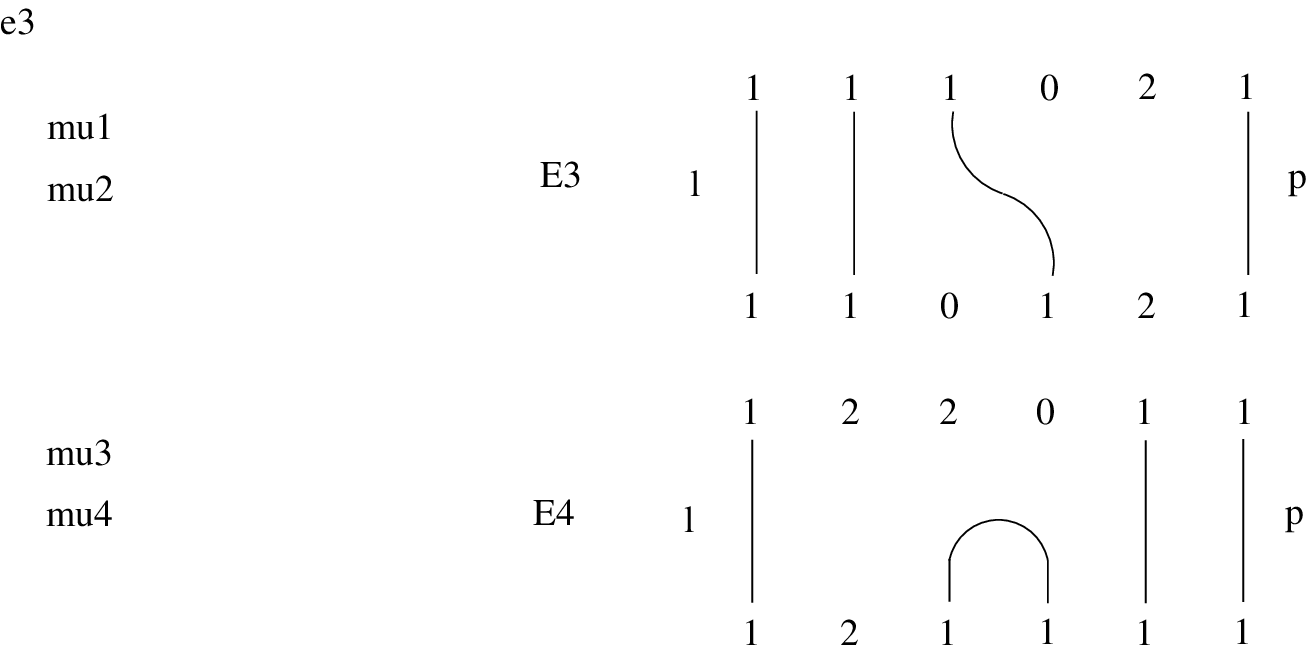}\caption{Examples of the functor $\mathcal{E}_i^{\mu}$.} \label{FunctorExample.figure}
\end{center}
\end{figure}

The functors $\mathcal{E}_i$, $\mathcal{F}_i$, and $\mathcal{K}_i$
defined by Khovanov and Huerfano naturally extend to our category
$\mathcal{C}$. Recall that $\mathcal{E}_i:\mathcal{C}\rightarrow
\mathcal{C}$ is defined to be the sum over all admissible $\mu$ of
the functors $\mathcal{E}_i^{\mu}: \mathcal{C}_{\mu}\rightarrow
\mathcal{C}_{\mu+\epsilon_i}$. If $\mu+\epsilon_i$ is not
admissible $\mathcal{E}_i^{\mu}$ is the zero functor. Otherwise,
define $\mathcal{E}_i^{\mu}$ to be tensoring with the
$(A_{\mu+\epsilon_i},A_{\mu})$-bimodule $\cF(T_{i}^{\mu})$ where
$T_{i}^{\mu}$ is the simplest flat tangle with bottom end points
corresponding to $\mu$ and top end points corresponding to
$\mu+\epsilon_i$. Figure~\ref{FunctorExample.figure} shows an
example of the functor $\mathcal{E}_i^{\mu}$.
The definition of $\mathcal{F}_i$ is similar. See
\cite{HK:LevelTwo} for details. Define $\mathcal{K}_i$ to be the
functor which shifts the grading of $M\in \mathcal{C}_{\mu}$ up by
$\mu_i-\mu_{i+1}$
$$\mathcal{K}_i(M)\define M\{\mu_i-\mu_{i+1}\}.$$

\begin{prop} \label{prop1} There are functor isomorphisms
\begin{equation}
 \label{functor-isom}
 \begin{array}{l}
 \cK_{i} \cK^{-1}_{i} \cong \mathrm{Id}
 \cong  \cK^{-1}_{i} \cK_{i}, \\
 \cK_{i} \cK_{j} \cong \cK_{j} \cK_{i}, \\
 \cK_{i} \cE_{j} \cong
 \cE_{j}\cK_{i}\{c_{i,j}\}, \\
 \cK_{i} \cF_{j} \cong
 \cF_{j} \cK_{i}\{-c_{i,j}\}, \\
 \cE_{i} \cF_{j} \cong  \cF_{j}\cE_{i}
 \hsp \mbox{ if } \hsp i \not= j, \\
 \cE_{i} \cE_{j} \cong \cE_{j}\cE_{i}
 \hsp\mbox{ if }\hsp |i-j|> 1, \\
 \cF_{i} \cF_{j} \cong \cF_{j}\cF_{i}
 \hsp\mbox{ if }\hsp |i-j|> 1, \\
 \cE_{i}^2 \cE_{j} \oplus
 \cE_{j} \cE_{i}^2 \cong
 \cE_{i} \cE_{j} \cE_{i}\{ 1\}  \oplus\cE_{i} \cE_{j} \cE_{i}\{ -1\}
 \hsp \mbox{ if } \hsp
  j = i \pm 1, \\
 \cF_{i}^2 \cF_{j} \oplus \cF_{j} \cF_{i}^2 \cong  \cF_{i}
 \cF_{j} \cF_{i} \{ 1\} \oplus\cF_{i} \cF_{j} \cF_{i} \{ -1\}
 \hsp \mbox{ if } \hsp
  j = i \pm 1,
\end{array}
\end{equation}
where

$c_{i,j} = {\left\{ \begin{array}{ll} 2 & \mathrm{ if }\hsp
  j = i, \\
 -1 & \mathrm{ if } \hsp j = i\pm 1, \\
   0 & \mathrm{ if} \hsp |j-i|>1. \end{array} \right. }$
\end{prop}

\begin{prop}
\label{prop2} For any  admissible $\mu$ there is an isomorphism of
functors in the category $\cC_\mu$
\begin{equation}
\label{more-fn}
 \begin{array}{l}
 \cE_i \cF_i \cong \cF_i \cE_i \oplus \id\{1\}\oplus \id \{ -1\}
\hspace{0.1in} \mbox{ if } \hspace{0.1in} (\mu_i,\mu_{i+1})=(2,0), \\
 \cE_i \cF_i \cong \cF_i \cE_i \oplus \id
\hspace{0.1in} \mbox{ if } \hspace{0.1in} \mu_i - \mu_{i+1} =1, \\
 \cE_i \cF_i \cong \cF_i \cE_i
\hspace{0.1in} \mbox{ if } \hspace{0.1in} \mu_i = \mu_{i+1}, \\
 \cE_i \cF_i \oplus \id \cong \cF_i \cE_i
\hspace{0.1in} \mbox{ if } \hspace{0.1in} \mu_i - \mu_{i+1}=-1, \\
 \cE_i \cF_i \oplus \id\{1\}\oplus \id \{ -1\} \cong \cF_i \cE_i
\hspace{0.1in} \mbox{ if } \hspace{0.1in} (\mu_i,\mu_{i+1})=(0,2).
 \end{array}
\end{equation}
\end{prop}

\begin{prop}
The functor $\cE_i$ is left adjoint to $\cF_i \cK_i^{-1}\{1\}$,
the functor $\cF_i$ is left adjoint to $\cE_i \cK_i\{1\}$, and
$\cK_i$ is left adjoint to $\cK_i^{-1}$.
\end{prop}

The above three propositions are from \cite{HK:LevelTwo}. They
work in our case without any modifications since the actions
happen away from the platform.

The Grothendieck group of $\cC$ is a $\mathbb{Z}[q,q^{-1}]$-module
where grading shifts correspond to multiplication by $q$. The
functors $\cE_i$, $\cF_i$, and $\cK_i$ are exact and commute with
grading shift action $\{1\}$. Exactness follows from left and
right projectivity of bimodule $\cF(T)$ for flat tangle $T$ in
section $2$. On the Grothendieck group level $\cE_i$, $\cF_i$, and
$\cK_i$ descend to $\mathbb{Z}[q,q^{-1}]$-linear endomorphisms
$[\cE_i]$, $[\cF_i]$, and $[\cK_i]$ respectively. Functor
isomorphisms in proposition \ref{prop1} and proposition
\ref{prop2} correspond to the quantum group relation (\ref{q-rel})
in $K(\cC)$. So we can view $K(\cC)$ as an $U_q(sl_N)$ module. It
follows from the bijection (\ref{bijection.equation}) that
$K(\cC)$ is isomorphic to $V$ as an $U_q(sl_N)$ module:

\begin{prop}
The Grothendieck group of $\cC$ is isomorphic to the irreducible
representation of $U_q(sl_N)$ with the highest weight
$\omega_k+\omega_{k+s}$
$$K(\cC)\otimes_{\mathbb{Z}[q,q^{-1}]}\mathbb{C}\cong V.$$
\end{prop}

\noindent Y.~Chen,  Department of Mathematics, University of
California, Berkeley, CA 94720.

\noindent yfchen@math.berkeley.edu

\end{document}